\documentclass[11pt]{amsart}

\usepackage{graphicx}
\usepackage{amsmath}
\usepackage{amssymb}
\usepackage{euscript}
\usepackage{caption}
\usepackage{subcaption}
\newtheorem{thm}{Theorem}[section]
\newtheorem{cor}[thm]{Corollary}
\newtheorem{lem}[thm]{Lemma}

\newtheorem{conj}[thm]{Conjecture}
\newtheorem{prop}[thm]{Proposition}
\theoremstyle{definition}

\theoremstyle{remark}

\numberwithin{equation}{section}

\newcommand{\DEF}[1]{\emph{#1}}

\begin{document}

\def \today {\number \day \ \ifcase \month \or January\or February\or
  March\or April\or May\or June\or July\or August\or
  September\or October\or November\or December\fi\
  \number \year}

\title{Egerv\'{a}ry Graphs: Deming Decompositions and Independence Structure}

\author{P. Mark Kayll}
\address{Department of Mathematical Sciences\\University
  of Montana\\Missoula MT 59812, USA}
\email{mark.kayll@umontana.edu}

\author{Craig E. Larson}
\address{Department of Mathematics and Applied
Mathematics\\Virginia Commonwealth University\\Richmond VA 23284, USA}
\email{clarson@vcu.edu}

\thanks{Partially supported by grants from the Simons Foundation
 (\#279367 to Mark Kayll and \#426267 to Craig Larson)}


\date{\today}

\maketitle

 \begin{abstract}

\noindent
We leverage an algorithm of 
Deming~[R.W. Deming, Independence numbers of graphs---an extension of
the Koenig-Egervary theorem, \textit{Discrete Math.}, \textbf{27}
(1979), no.\,1,  23--33; MR534950]  to
decompose a matchable graph into subgraphs with a
precise structure: they are either spanning even subdivisions of
blossom pairs, spanning even subdivisions of the complete graph $K_4$,
or a K\H{o}nig-Egerv\'{a}ry graph. In each case, the 
subgraphs have perfect matchings; in the first two cases, their
independence numbers are one less than their matching numbers, while
the independence number of the KE subgraph equals
its matching number.  This decomposition refines previous results
about the independence structure of an arbitrary graph and leads to
new results about $\alpha$-critical graphs.

\medskip
\noindent
\textbf{Keywords:}  matching, independence, K\H{o}nig-Egerv\'{a}ry,
Egerv\'{a}ry, Birkhoff-von Neumann 

\medskip
\noindent
MSC2020:\ Primary
05C70,   
05C69;  
Secondary
05C75,   
05C85.   

   \end{abstract}

%
\renewcommand{\thefootnote}{}
\footnotetext{\copyright~2022 by the authors}

\section{Introduction}

K\H{o}nig-Egerv\'{a}ry graphs generalize bipartite graphs: they are
defined by the condition that their independence number $\alpha$
and their matching number $\nu$ sum to their order $n$.
These graphs have the property that they can be
identified efficiently and their independence numbers can be
calculated efficiently. This investigation is motivated by an attempt
to generalize K\H{o}nig-Egerv\'{a}ry graphs to a larger class of
graphs with similar attractive properties.  

In this paper, we leverage an algorithm of Deming~\cite{Demi79} to
produce a useful 
decomposition of a matchable graph into subgraphs with a
precise structure: they are either spanning even subdivisions of
blossom pairs, spanning even subdivisions of the complete graph $K_4$,
or a K\H{o}nig-Egerv\'{a}ry graph. In all cases, the 
subgraphs have perfect matchings; in the first two cases, their
independence numbers are one less than their matching numbers, while
the independence number of the K\H{o}nig-Egerv\'{a}ry subgraph equals
its matching number.  This decomposition refines previous results
about the independence structure of an arbitrary graph and leads to
new results about $\alpha$-critical graphs. 

\subsection{Matching Structure, K\H{o}nig-Egerv\'{a}ry \& Egerv\'{a}ry Graphs}
\label{sec:matching}

An \textit{independent set} is a set of vertices that is pairwise
nonadjacent; the \textit{independence number} $\alpha$ of a graph is
the cardinality of a largest 
independent set. A \textit{matching} is a set of pairwise non-incident
(or independent) edges; the \textit{matching number} $\nu$ is the
cardinality of a largest matching. The \textit{order} $n$ of a graph
is the number of its vertices. A \textit{K\H{o}nig-Egerv\'{a}ry} graph
(henceforth, \textit{KE} graph) is defined by the 
condition that $\alpha+\nu=n$. While the definition of an Egerv\'{a}ry
graph and our initial theorems apply to graphs with perfect matchings, we
will ultimately show that many of our results extend to arbitrary graphs.  

Egerv\'{a}ry graphs generalize KE graphs.
A graph is \textit{Egerv\'{a}ry} if it has a perfect matching but it
doesn't contain a spanning subgraph whose components are 
independent edges together with at least one odd cycle (parity 
considerations show that ``at least one'' really means 
``a positive even number of''). These graphs are also known
as `non-Edmonds' graphs (e.g., in \cite{Kayl10}) and `Birkhoff-von
Neumann' (BvN) graphs (in \cite{CarvKothWangLin20}). The definition
above is the most convenient one for our purposes, but it's worth
noting that these graphs also admit a polyhedral description depending
on their perfect matching polytope. See \cite{Kayl10} for details.
That Egerv\'{a}ry graphs indeed generalize KE graphs will become
apparent with the statement of Corollary~\ref{cor:kayll-KE-implies-eger}.

One should observe the similarity of our definition of
(non-)Egerv\'{a}ry graphs with the notion of a 
\textit{basic perfect $2$-matching} of a graph $G$, 
namely a spanning subgraph whose components are 
independent edges and odd cycles; see, e.g., \cite{LovaPlum86}. For
matchable graphs, being Egerv\'{a}ry is almost the negation of
admitting a basic perfect $2$-matching, the difference stemming from
the stipulation `at least one' in the first definition.
We don't make much use of this connection, but it does appear in the
proof of Lemma~\ref{lem:Andrasfai-helper} below.

Egerv\'{a}ry graphs are also related to \textit{tangled} graphs, i.e., 
those graphs containing no two vertex-disjoint odd cycles; see, e.g., 
\cite{Slil07} or \cite{KawaOzek13} and the references cited there.
The definitions show that matchable tangled graphs are Egerv\'{a}ry,
but the converse is not true. In fact, tangled graphs admit a
characterization in terms of `super-Egerv\'{a}ry' graphs,
whose definition is beyond our scope. But we note that the omitted 
definition generalizes the polyhedral description of Egerv\'{a}ry
graphs to which we alluded above. We hope to return to this thread in 
a follow-up article.  

\begin{center}
\begin{figure}
\includegraphics[width=9cm]{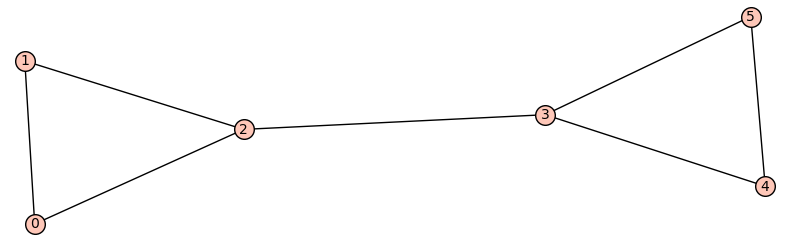}
\caption{The graph $T$, which is not KE ($\alpha=2$, $\nu=3$, and
  $n=6$); it is also not Egerv\'{a}ry (the two disjoint $3$-cycles span $T$).} 
\label{T}
\end{figure}
\end{center}

The class of Egerv\'{a}ry graphs includes, for instance, the complete graph $K_4$
but not the graph $T$ consisting of two $3$-cycles and a
single edge incident to one vertex of each (see Figure~\ref{T}).

Deming~\cite{Demi79} and Sterboul~\cite{Ster79} (independently) gave
the first characterizations of KE graphs. In particular, they showed
that this property is in co-NP. There are now many characterizations of
KE graphs; the following one is particularly useful for the current
investigation. Following \cite{LovaPlum86}, we call a subgraph $H$ of
a matchable graph $G$ \textit{nice} if $G-V(H)$ has a perfect matching. 

\begin{thm}[Lov\'asz \& Plummer, \cite{LovaPlum86}]
\label{KE-char-LP86}
A graph with a perfect matching is KE  if and only if it does not
contain an even subdivision of either $K_4$ or $T$ as a nice subgraph. 
\end{thm}

\noindent
Even subdivisions of $T$ are sometimes called \textit{blossom pairs}. 
It's easy to check that  no blossom pair is KE (and that neither is
any even $K_4$-subdivision).
Theorem~\ref{KE-char-LP86} shows that for matchable graphs, these are
the only obstructions.

The proof of the next theorem imitates the proof of
Theorem~\ref{KE-char-LP86}; since it 
restricts the KE characterization to fewer forbidden subgraphs, it
generalizes the characterized class of graphs. 

\begin{thm}
A graph with a perfect matching can be covered by a collection of
independent edges and a positive number of odd cycles  if and only if
it contains an even subdivision of $T$ as a nice subgraph. 
\end{thm}

\begin{proof}
Let $G$ be a graph with a perfect matching $M_1$. If $G$ contains an 
even subdivision of $T$ as a nice subgraph, then $G$ can certainly 
be covered by a collection of independent edges and a positive 
number of odd cycles.

Now suppose that $G$ can be covered by a collection $M_0$ of
independent edges and a positive number of odd cycles $C_1,\ldots,C_k$.
Then $M_0\cup M_1$ consists of the edges in $M_0\cap M_1$, some
alternating (even) cycles, and some paths connecting vertices in the
odd cycles 
$C_1,\ldots,C_k$. These paths start and end with edges of $M_1$.
At least one such path must connect vertices in different cycles since
they cannot pair off the vertices within a single odd cycle.  

So let $P$ be a path-component of $M_0\cup M_1$ connecting two odd
cycles $C_i$ and $C_j$. The subgraph $P\cup C_i\cup C_j$ is an even
subdivision of $T$, and it is also a nice subgraph of $G$. 
\end{proof}

\begin{cor}
\label{Egervary_char-even-T}
A graph with a perfect matching is Egerv\'{a}ry 
 if and only if it does not contain an even subdivision of $T$ as a
 nice subgraph. 
\end{cor}

\begin{center}
\begin{figure}
\includegraphics[width=5cm]{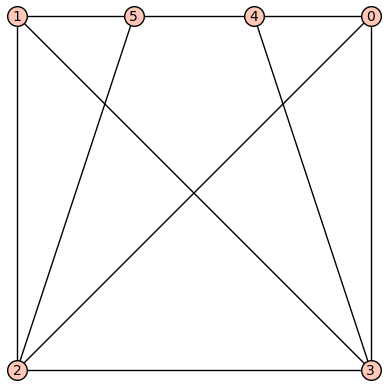}
\caption{This matchable graph contains a spanning even subdivision of
  $K_4$ (so is not KE by Theorem~\ref{KE-char-LP86}) and a nice even
  subdivision of $T$ (so is not Egerv\'{a}ry by Corollary~\ref{Egervary_char-even-T}).}  
\label{k4_problem}
\end{figure}
\end{center}

\noindent
Figure~\ref{k4_problem} illustrates Theorem~\ref{KE-char-LP86} and
Corollary~\ref{Egervary_char-even-T}.
These results also give a new proof of the following
theorem, originally proved by the first author of the present work.

\begin{cor}[\cite{Kayl10}] 
\label{cor:kayll-KE-implies-eger}
If a graph with a perfect matching is KE, then it is Egerv\'{a}ry. 
\end{cor}

\section{Deming Subgraphs}
\label{sec:deming}

We shall show that every graph with a perfect matching admits a
decomposition into subgraphs with attractive properties. Given such a
graph $G$, Deming's Algorithm~\cite{Demi79} either
produces a maximum independent set $I$ that certifies $G$
being KE or produces an even subdivision of $K_4$ or $T$ as a
nice subgraph. We can push this algorithm further so that, when
$G$ is not KE, we get either a `Deming-$K_4$ subgraph'---with 
a spanning even $K_4$-subdivision---or a
`Deming-BP subgraph' with a spanning even $T$-subdivision, 
both of which are, in a precise sense, \textit{almost}-KE. We describe
these first and then discuss our extension of Deming's Algorithm. 

A graph $K$ is a \textit{Deming-$K_4$ graph} if it contains a
spanning even subdivision of $K_4$ and $K-\{x,y\}$ is KE for every
edge $xy$ in \textit{some} perfect matching of $K$. 
Analogously, a graph $B$  is a \textit{Deming-BP graph} 
(for `Deming-Blossom-Pair') if it
contains a spanning even subdivision of  the graph $T$ (in Figure~\ref{T})
and $B-\{x,y\}$ is KE for every edge $xy$ in some perfect matching of $B$.  
Finally, a graph $D$ is a \textit{Deming graph} if it is either a
Deming-$K_4$ graph or a Deming-BP graph.

\begin{center}
\begin{figure}
\includegraphics[width=10cm]{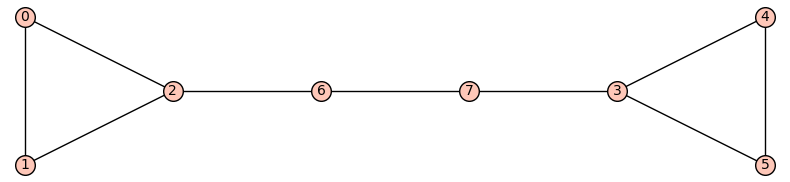}
\caption{A Deming-BP graph $B$: 
  an even subdivision of $T$ such that, for every edge $xy$ of $B$ in
  some perfect matching, the vertex-deleted graph $B-\{x,y\}$  is KE.} 
\label{T_sub}
\end{figure}
\end{center}

Figure~\ref{T_sub} shows an example of a Deming-BP graph. Notice that for
the center edge 6--7, the subgraph $B-\{6,7\}$ is not KE, but 6--7 is
not in a perfect matching of $B$. It is also easy to check that the
Petersen graph $P$ is a Deming-BP graph: it has a pair of 5-cycles
connected by a single edge as a spanning subgraph, every edge $xy$ is
in some perfect matching, and $P-\{x,y\}$ is KE.  

Note that whether a graph $D$ with a perfect matching $M$ is Deming
can be efficiently determined: run Deming's Algorithm on $D$ (with
resect to $M$). Deming's Algorithm will produce either a maximum
independent set when $D$ is KE, or a nice even subdivision
of either $K_4$ or $T$ when it is not. If $D$ is not
KE and the produced even subdivision spans $D$, then the graph may be
Deming. It now remains to identify which edges $xy$ are in some
perfect matching (by checking, for instance if $D-\{x,y\}$ has a
perfect matching), and then, for each such edge $xy$, checking if
$D-\{x,y\}$ is KE. Of course, this algorithm leverages the fact that
it is efficient to find a maximum matching in a graph \cite{Edmo65}. 

One output of Deming's Algorithm is subgraphs which have a spanning
even subdivision of $K_4$. We dig  
deeper into the subgraphs produced by Deming's Algorithm. Even
subdivisions of $K_4$ have the property that every edge is in a
perfect matching and that the removal of any edge yields a KE graph
(and they are `$\alpha$-critical'). The subgraphs produced by
Deming's Algorithm may have `extra' edges. The definition here is
designed to efficiently reduce a given output of Deming's Algorithm
into a similar substructure with even nicer properties. 

Deming-$K_4$ graphs have even subdivisions of $K_4$; as these are
central in this paper, we record some relevant properties.  

\begin{prop}\label{perfect}
An even subdivision of $K_4$ has a perfect matching.
\end{prop}

\begin{proof}
An even subdivision $K$ of  $K_4$
consists of four `corner' vertices and six odd paths between
them. Let $V(K_4)=\{v_1,v_2,v_3,v_4\}$ and let the edges of $K_4$ be
$e_{i,j}$, one for each pair of distinct vertices $v_i$, $v_j$ in
$K_4$. Call the corresponding corner vertices in $K$:
$v_1$, $v_2$, $v_3$, $v_4$. And for each pair $v_i$, $v_j$ of distinct
corner vertices, let $P_{i,j}$ be the odd path from $v_i$ to $v_j$ formed by
subdividing $e_{i,j}$ in $K_4$ an even number of times. The odd-length
paths $P_{1,2}$ and $P_{3,4}$ can be matched as can the odd-length
paths strictly between vertices $v_1$ and $v_3$ on
$P_{1,3}$, strictly between $v_2$ and $v_4$ on
$P_{2,4}$, strictly between $v_1$ and $v_4$ on
$P_{1,4}$, and (finally) strictly between $v_2$ and $v_3$ on
$P_{2,3}$.
\end{proof}

\begin{cor}
\label{cor:all_K4edges_in_PM}
Every edge of an even subdivision of $K_4$ is in some perfect matching. 
\end{cor}

\begin{cor}\label{DK4perfect}
A Deming-$K_4$ graph has a perfect matching.
\end{cor}

\begin{cor}\label{Dperfect}
A Deming graph has a perfect matching.
\end{cor}

\begin{proof}
Deming-$K_4$ graphs have a perfect matchings by the preceding
corollary. A Deming-BP graphs has a perfect matching as it contains a
spanning even $T$ subdivision: the odd path connecting the blossom
tips has a perfect matching, and the remaining non-blossom-tip
vertices in both blossoms also form odd paths. 
\end{proof}

The following related results are needed for the discussion of
$\alpha$-critical graphs in Section~\ref{sec:structure}. 

\begin{prop}\label{G-e}
If $K$ is an even subdivision of $K_4$ and $e$ is an edge of $K$, then
$K-e$ is KE. 
\end{prop}

\begin{proof}
We again view $K$ as having corner vertices $v_1,v_2,v_3,v_4$ with
odd-length paths $\{P_{i,j}\}$ between them.  Assume that $e=xy$ is on
$P_{1,2}$ with $x$ closer to $v_1$ than is $y$. We shall first show
that $K$ contains a perfect matching not containing $e$. 

With $P_{1,2}$ being of odd length, the parities of the distances from
$x$ to $v_1$ and $y$ to $v_2$ in $P_{1,2}$ are necessarily the same.

If the parities are both odd, then taking the remaining edge of
$P_{1,2}$ incident to vertex $x$ and every alternate edge from $x$ to
$v_1$ will include an edge that saturates $v_1$. Include
these edges in a matching $M$. Similarly, taking the remaining edge of
$P_{1,2}$ incident to $y$ and every other edge from $y$ to $v_2$ will
include an edge that saturates $v_2$. Add the edge incident to $v_3$
to $M$ along with every other edge from $v_3$ to $v_4$; this will
include an edge that covers $v_4$. Now all corner vertices are covered
by $M$ as well as every vertex on $P_{1,2}$ and $P_{3,4}$. It remains
to add edges that saturate the interior vertices of the other four
paths $P_{1,3}$, $P_{2,4}$, $P_{1,4}$, and $P_{2,3}$. In each case,
the remaining paths between (and not including) the corresponding
corner vertices  have odd length, and these vertices can all be
covered by matching edges added to $M$. 

If the distance parities from $x$ to $v_1$ and $y$ to $v_2$ in
$P_{1,2}$ are both even, then taking the remaining edge of
$P_{1,2}$ incident to $x$ and every other edge from $x$ to $v_1$ will
\emph{not} include an edge that saturates $v_1$. Include these edges
in $M$. Similarly, taking the remaining edge of $P_{1,2}$ incident to
$y$ and every other edge from $y$ to $v_2$ will \emph{not} include an
edge that saturates $v_2$. Also include in $M$ a matching that covers
the odd-length paths $P_{1,4}$ and $P_{2,3}$. These matching edges
will include all corner vertices. It remains to add edges that
saturate the interior vertices of the other three paths $P_{1,3}$,
$P_{2,4}$, and $P_{3,4}$. In each case, the remaining paths between
(and not including) the corresponding corner vertices have odd
lengths, and these vertices can all be covered by matching edges  added to $M$.

So in both the odd- and even-parity cases, $K-e$ has a perfect
matching. Using $M$ to denote that matching in either case, 
we shall construct an independent set $I$ with $|I|=|M|$. As this
condition implies that $K-e$ is KE, arranging for it will suffice to
complete the proof. 

In the first case, $x$ is an odd distance to $v_1$.
Along $P_{1,2}$, build $I$ by including $x$, including every alternate vertex
along the $(x,v_1)$-segment of $P_{1,2}$, and  including $y$ and every alternate
vertex along the $(y,v_2)$-segment of $P_{1,2}$. Continuing around 
the cycle  $P_{1,2}$ followed by  $P_{2,3}$,  $P_{3,4}$, and
$P_{4,1}$, extend $I$ by including $v_2$'s neighbor along $P_{2,3}$
and then every alternate vertex (moving away from $v_2$) up to but not
including $v_1$. If that cycle is $\EuScript{C}$, then around
$\EuScript{C}$, we so far have selected one vertex for $I$ from each
$M$-edge of $\EuScript{C}$---including $v_3$ but not $v_1$, $v_2$, or
$v_4$. It remains to make analogous choices for vertices of $I$ along 
$P_{1,3}$ and $P_{2,4}$. Because $v_3\in I$, we extend $I$ on
$P_{1,3}$ by choosing all vertices at even distance from $v_3$ along
$P_{1,3}$; notice that this puts $v_1$'s neighbor on $P_{1,3}$ into
$I$---and not $v_1$. This is important because $v_1$'s neighbor on
$P_{1,2}$ is already in $I$, and we're building an independent set. 
Because $v_2,v_4\not\in I$, we have some choice for extending $I$ on
$P_{2,4}$. For definiteness, we select for $I$ all vertices at odd
distance from $v_2$ along $P_{2,4}$, except for $v_4$ (whose neighbors
along both $P_{4,1}$ and $P_{3,4}$ are already in $I$). One can check
that the described set $I$ is independent and contains exactly one
vertex of each $M$-edge when $M$ is the perfect matching defined for
this case. Thus,   $|I|=|M|$ as desired.

In the second case, $x$ is an even distance to $v_1$ and $y$ is an
even distance to $v_2$. Similarly to the odd-parity case above, we may
again construct an independent set $I$---necessarily containing $x$
and $y$---with $|I|=|M|$, for $M$ now the perfect matching defined for
the even case. Here we omit a detailed description.
\end{proof}

One motivation for our definition of Deming-$K_4$ and
Deming-BP graphs is that an arbitrary graph can be
efficiently decomposed into subgraphs of three types with a precise
independence and matching structure. In particular,  Deming graphs
satisfy $\alpha=\nu-1$. 

\begin{thm}\label{thm:deming}
If $D$ is a Deming graph, then $\alpha(D)=\nu(D)-1$.
\end{thm}

\begin{proof}
Corollary~\ref{DK4perfect} shows that $D$ has a perfect matching,
and thus $\nu(D)=n(D)/2$. 
Deming graphs are not KE, and thus $\alpha(D) < n(D)/2=\nu(D)$; and so 
$\alpha(D)\leq\nu(D)-1$. Let $xy$ be an edge in a perfect matching of
$D$. Consider the graph $D'=D-\{x,y\}$ (formed by removing vertices $x$
and $y$ and all incident edges). Since $xy$ is in a perfect matching
of $D$, the subgraph $D'$ also has a perfect matching, and
$\nu(D')=\nu(D)-1$. By definition, $D'$ is KE; this and
the fact that $D'$ is matchable give
$\alpha(D')=\nu(D')$. Then
$\alpha(D)\geq\alpha(D')=\nu(D')=\nu(D)-1$. 
So $\alpha(D)=\nu(D)-1$.
\end{proof}

\begin{center}
\begin{figure}
\includegraphics[width=5cm]{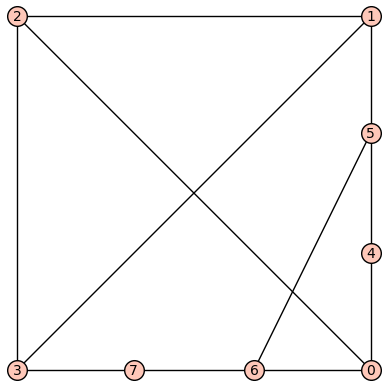}
\caption{a Deming-$K_4$ graph $K$ that is not an even subdivision 
of $K_4$. Notice that neither $K-\{\text{5--6}\}$ nor
$K-\{\text{0--6}\}$ is KE because both contain 
spanning even $K_4$-subdivisions. Theorem~\ref{thm:evensub}
characterizes when all such edge-deleted subgraphs are KE.}
\label{k4notegervary}
\end{figure}
\end{center}

By definition, a Deming-$K_4$ graph $K$ may contain edges not participating
in an even $K_4$-subdivision spanning it. For example, in
Figure~\ref{k4notegervary}, the edge 5--6 lies outside the 
$K_4$-subdivision with corner vertices $\{0,1,2,3\}$. So deleting this
edge produces a non-KE graph. The same is true for the edge-deleted
subgraph $K-\{\text{0--6}\}$, which can be seen to contain a different
spanning even $K_4$-subdivision. The next result provides
a way to test whether $K$ contains any of these superfluous edges,
namely by testing whether each of $K$'s single-edge-deleted subgraphs
is KE.

\begin{thm}\label{thm:evensub}
A Deming-$K_4$ graph $K$ is an even subdivision of $K_4$  if and only
if for every edge $e$ of $K$, the graph $K-e$ is KE. 
\end{thm}

\begin{proof}
Denote a perfect matching in $K$ by $M$.
First suppose that $K$ is an even subdivision of $K_4$,
with corner vertices $v_1$, $v_2$, $v_3$, $v_4$   
and paths $\{P_{i,j}\}$ defined as in the proof of Proposition~\ref{perfect}.
Let $e$ be an edge of $K$. We showed in Proposition~\ref{G-e} that
$K-e$ is KE.   

Next suppose that $K-e$ is KE for every edge $e$ of $K$. 
By the definition of a Deming-$K_4$ graph, $K$  has a spanning even
$K_4$-subdivision. We'll show that the subdivision edges are in
fact the only edges of $K$. Again, let the corner vertices of this
subdivision be $\{v_1,v_2,v_3,v_4\}$ and the six odd-length
paths between them be $\{P_{i,j}\}$.
Suppose that $xy$ is an edge not in any of the paths $P_{i,j}$. We may
assume that $x\in P_{1,2}$.  

As a first case, suppose that $x=v_1$ and $y \in P_{1,2}$. 
Either the distance from $x$ to $y$ on $P_{1,2}$ is odd or it is
even. If this distance is even, then there is a nice even subdivision of
$T$ in $K$. There are odd cycles: (1) the path from $x$ to $y$ in $P_{1,2}$
and back to $x$; and (2) the path $P_{2,3}$, followed by $P_{3,4}$
and then $P_{4,2}$, with blossom tips $x=v_1$ and $v_3$ joined by the
odd path $P_{1,3}$. When the distance from $x$ to $y$ on
$P_{1,2}$ is odd, at least one of the edges $e$ on the path from $x$ to
$y$ in $P_{1,2}$ is not in $M$. Now $K-e$ still contains a nice even
subdivision of $K_4$ and thus cannot be KE. 

Next we assume that while both $x$ and $y$ are in
$P_{1,2}$, neither of them is $v_1$ or $v_2$; assume that $x$ is closer to
$v_1$ than $y$. Again, either the distance from $x$ to $y$ on $P_{1,2}$
is odd or it is even. If this distance is even, then there is a nice
even subdivision of $T$ in $K$. Assume the distance from $x$ to $v_1$
along $P_{1,2}$ is even (either it is or the distance from $y$ to
$v_2$ must be). There are odd cycles: (1) the path from $x$ to $y$ in
$P_{1,2}$ and back to $x$; and (2) the path $P_{2,3}$, followed by
$P_{3,4}$ and then $P_{4,2}$, with blossom tips $x$ and $v_3$ joined
by the odd path from $x$ to $v_1$, finally followed by $P_{1,3}$. If
the distance from $x$ to $y$ on $P_{1,2}$ is odd, then at least one of
the edges $e$ on the path from $x$ to $y$ in $P_{1,2}$ is not in
$M$. Now $K-e$ still contains a nice even subdivision of $K_4$ and
thus cannot be KE. 

The next case we consider is when $x$ and $y$ are on adjacent
subdivision paths, say $P_{1,2}$ and $P_{2,3}$. Either the distance
from $x$ to $v_1$ in $P_{1,2}$ and from $y$ to $v_3$ in $P_{2,3}$ have
the same parity, or they do not. If the parities are the same, then there is a
nice even subdivision of $T$ in $K$. There are odd cycles: (1) the
path from $x$ to $v_2$ in $P_{1,2}$, from $v_2$ to $y$ in $P_{2,3}$, and
back to $x$; and (2) the path $P_{1,3}$, followed by $P_{3,4}$, and
then $P_{4,1}$, with blossom tips $v_2$ and $v_4$ joined by the odd
path $P_{2,4}$. In the special case where $x=v_1$ and $y=v_3$ (so the
distances are 0 and the parities are even), at least one of the edges
$e$ on the path from $v_1$ to $v_3$ in $P_{1,3}$ is not in $M$. Now
$K-e$ still contains a nice even subdivision of $K_4$ and thus cannot
be KE. If the parities from $x$ to $v_1$ in $P_{1,2}$ and from $y$ to
$v_3$ in $P_{2,3}$ are different (assume the distance from $y$ to
$v_3$ in $P_{2,3}$ is odd), at least one of the edges $e$ on the path
from $x$ to $v_2$ in $P_{1,2}$ and from $v_2$ to $y$ in $P_{2,3}$ is
not in $M$. Now $K-e$ still contains a nice even subdivision of $K_4$
and thus cannot be KE. 

The last case we consider is when $x$ and $y$ are on opposite
subdivision paths, say $P_{1,2}$ and $P_{3,4}$.  Let $e$ be an edge in
$P_{2,3}$. By assumption, $K-e$ is KE. We will show that in
fact (and contrary to our assumption), $K$ contains a nice even
subdivision of $K_4$. We may assume that neither $x$ nor $y$ is a
corner vertex---otherwise these fall into cases that we handled
previously. The distances from $x$ to the corner points of $P_{1,2}$
necessarily have different parities. Similarly, the distances from $y$
to the corner points of $P_{3,4}$ have different
parities. Assume that the distances from $x$ to $v_1$ and from $y$ to
$v_4$ are both odd. We now construct an even subdivision of $K_4$
with corner points $v_1$, $x$, $y$, and $v_4$. There is an odd path
from $x$ to $v_4$ consisting of the even path from $x$ to $v_2$ in
$P_{1,2}$ followed by $P_{2,4}$. And there is also an odd path from
$y$ to $v_3$ consisting of the even path from $y$ to $v_3$ in
$P_{3,4}$ followed by $P_{3,1}$.  
\end{proof}

A graph $G$ is \textit{$\alpha$-critical} if
$\alpha(G-xy)=\alpha(G)+1$ for each edge $xy$ in $G$. 
(We include $K_1$ as 
(vacuously) $\alpha$-critical; it is not only
true but also makes the statements of some of our results cleaner.)
The study of $\alpha$-critical graphs has a long history
\cite{ErdoGall61,BeinHaraPlum67,Lova79a,SeweTrot93}. We record one
fact here and return to this topic in Section~\ref{sec:structure}. 

\begin{cor}\label{cor:alphacriticalk4}
A Deming-$K_4$ graph $K$ is an even subdivision of $K_4$  if and only
if $K$ is $\alpha$-critical.  
\end{cor}

\begin{proof}
Suppose that $K$ is an even subdivision of $K_4$. It immediately
follows from Theorem~\ref{thm:evensub} that $K$ is $\alpha$-critical. 

Suppose then that $K$ is a Deming-$K_4$ graph and is
$\alpha$-critical. Then $K$ has a spanning
even $K_4$-subdivision $F$. Note that $F$ is also a Deming-$K_4$
graph and that $K$ and $F$ must both have the same matching
number. Theorem~\ref{thm:deming} implies both that
$\alpha(K)=\nu(K)-1$ and that $\alpha(F)=\nu(F)-1$. So
$\alpha(K)=\alpha(F)$. Now suppose that $xy$ is a non-$F$ edge of
$K$. Clearly 
$\alpha(K)\leq \alpha(K-xy)\leq \alpha(F-xy)=\alpha(F)=\alpha(K)$.
So we have both that $\alpha(K)=\alpha(K-xy)$ and that
$\alpha(K-xy)=\alpha(K)+1$. Since this is impossible, $K$ cannot
contain a non-$F$ edge. 
\end{proof}

\begin{prop}
\label{Ev-T-subdiv-have-uniq-PM}
If $K$ is an even $T$-subdivision, then it has a unique perfect
matching $M$, and for every edge $e$ not in $M$, the graph $K-e$ is KE
and  has a unique perfect matching. 
\end{prop}

\begin{proof}
Let $K$ be formed from two odd cycles $C_1$, $C_2$ with respective
blossom tips $v$, $w$ joined by an odd-length path $P$. Let
$V(C_1)=\{v=v_1,\ldots,v_{2k+1}\}$ and
$V(C_2)=\{w=w_1,\ldots,w_{2l+1}\}$. Since every perfect matching $M$ of
$K$ must include the edges of $P$ that cover $v$ and $w$ (and every
alternate edge between them), the remaining uncovered vertices of $C_1$
form the odd path from $v_2$ to $v_{2k+1}$. As there is a unique
matching $M_1$ covering these vertices, $M$ must include the edges
of $M_1$. Similarly, the remaining
uncovered vertices in $C_2$ form the odd path from $w_2$ to
$w_{2l+1}$, and there is a unique matching covering these
vertices. Thus, $M$ must include these edges as well. So $G$ has a
unique perfect matching. 

Now let $e=xy$ be a non-$M$ edge. There are two cases: either $e$ is
on the path $P$ or $e$ belongs to one of the
cycles (which we may assume to be $C_1$). Suppose first that $e$
belongs to $P$. Assume that $x$ is the end of $e$
closer to $v$ (so $y$ is the end closer to $w$). So there is an
odd-length path from $v$ to $x$ and one from $y$ to $w$. The graph
$K-e$ then comprises two components, each consisting of an odd cycle 
together with an odd-length path adjoined to the blossom tip. Since
both of these components are KE, so too is $K-e$. 

Now suppose that $e$ is a non-$M$ edge in $C_1$. Assume that $x$
precedes $y$ on the path from $v_2$ to $v_{2k+1}$. Then $x$ must equal
$v_i$ for some odd index $i$ (and thus $y=v_{i+1}$ with $i+1$
even). Let $I$ be the independent set containing the vertices
$w_2,w_4,\ldots,w_{2l}$, the vertices at odd distance from $w$ on $P$
(including $v$), the vertices at even distance from $v$ on the path
$v=v_1,v_2,\ldots,x$ (including $x$), and finally the vertices at even
distance from $y$ along the path $y=v_{i+1},\ldots v_{2k}$ (including
$y$). The set $I$ thus includes exactly one end of each $M$-edge;
whence $|I|=|M|$ and $G-e$ is KE.

Of course, $K-e$  has a unique perfect matching because $K$ does. 
\end{proof}

\begin{cor}\label{cor:alphacriticaldeming}
If $D$ is an $\alpha$-critical  Deming graph, then $D$ is an even
subdivision of $K_4$. 
\end{cor}

\begin{proof}
First note that the graph $T$ is not $\alpha$-critical; neither is
\textit{any} blossom pair $D$: if $x$ and $y$ are blossom tips and
$xx'$ is the edge incident to $x$ on the path from $x$ to $y$, then
$\alpha(D)=\alpha(D-xx')$. Extending this observation, we can show
that no Deming-BP graph $D$ is $\alpha$-critical. For such a $D$ is
either a blossom pair (in which case we're done) or 
contains an edge $uv$ that is not an edge in $D$'s spanning blossom
pair. Let $D'=D-uv$. It is easy to check that $D'$ is also a Deming-BP graph.  
Theorem~\ref{thm:deming} then implies that $\alpha(D)=\nu(D)-1$ and
$\alpha(D')=\nu(D')-1$. Since $D'$ has a perfect matching, whence
$\nu(D)=\nu(D')$, it follows that $\alpha(D)=\alpha(D')=\alpha(D-uv)$,
and thus indeed $D$ is not $\alpha$-critical.

So if $D$ is an $\alpha$-critical Deming graph, it must be a
Deming-$K_4$ graph, and now Corollary~\ref{cor:alphacriticalk4} shows
that $D$ is an even subdivision of $K_4$.
\end{proof}

\section{Extending Deming's Algorithm and Deming Decompositions}
\label{sec:extending}

Given a graph $G$ with a perfect matching, Deming's original
algorithm~\cite{Demi79} 
either produces a maximum independent set $I$ certifying
that $G$ is KE or produces an even subdivision of either $K_4$ or $T$
as a nice subgraph. In the later cases, we test if the produced
obstructions are Deming graphs and, if not, re-apply Deming's
Algorithm to certain subgraphs. 

Suppose that Deming's Algorithm is applied to a graph $G$, and the
output is an obstruction $H$, i.e., an even subdivision of $T$ or
$K_4$. In either case, we find the edges $xy$ which are in some
perfect matching of $H$ and, for each of these, test whether $H-\{x,y\}$ is
KE by re-applying Deming's Algorithm. If there is an edge $xy$ such
that $H-\{x,y\}$ is not KE, then this graph must itself contain an
obstruction $H'$: a nice even subdivision of $T$ or $K_4$ . We then
iterate on $H-\{x,y\}$. Since we remove two vertices on each
iteration, we must terminate (in no more than 
$n(H)/2$ steps). If $H$ is not KE, then this final output obstruction
$D$ must be a Deming graph. 

Now the vertices of $D$ can be removed from $H$. The resulting graph
$H-D=G[V(H)\setminus V(D)]$ must itself have a perfect matching (as
the obstruction $D$ contains a nice even subdivision of  either $T$ or
$K_4$), and this \textit{extended Deming Algorithm} can be applied to
the subgraph $H-D$. Since we remove at least four vertices on each
iteration, we must terminate (in no more than $n(H)/4$
steps). The final output graph, possibly 
null, is necessarily KE. So
this process results in a finite sequence of Deming graphs
together with a single KE graph whose vertices partition the vertices
of the parent graph $G$ (allowing some parts to be empty). 

\begin{figure}
\begin{subfigure}{0.35\textwidth}
\includegraphics[width=0.9\linewidth]{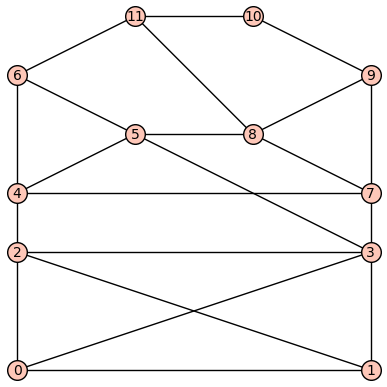}
\caption{$G$}
\label{sub1}
\end{subfigure}
\begin{tabular}{ccc}
\begin{subfigure}{0.35\textwidth}
\includegraphics[width=0.9\linewidth]{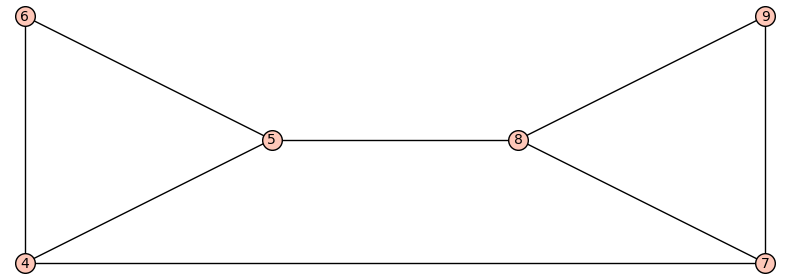}
\caption{$B$}
\end{subfigure}
&
\begin{subfigure}{0.25\textwidth}
\includegraphics[width=0.9\linewidth]{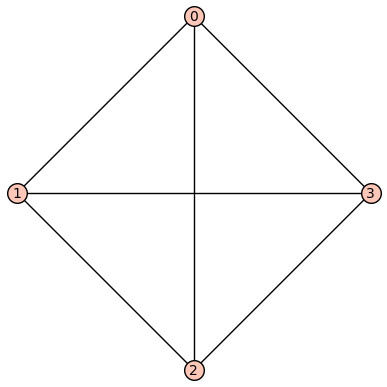}
\caption{$K$}
\end{subfigure}
&
\begin{subfigure}{0.15\textwidth}
\includegraphics[width=0.9\linewidth]{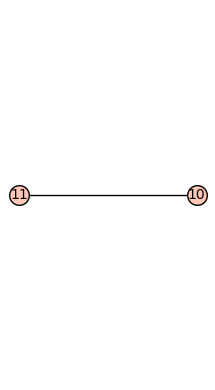}
\caption{$R$}
\end{subfigure}
\end{tabular}
\caption{a Deming decomposition of $G$ in (\textsc{a}) is $\{\{B\},\{K\},R\}$, where $B$ is a Deming-BP subgraph, $K$ is a Deming-$K_4$ subgraph, and $R$ is KE}
\label{Deming-decomp-examp}
\end{figure}

It is easy either to keep track, or efficient to test directly,
whether the output Deming graphs are Deming-$K_4$ subgraphs or
Deming-BP subgraphs, all of which defines a decomposition of $G$ into
subgraphs with well-defined properties. 
We call this a `Deming decomposition'; to formalize, a 
\textit{Deming decomposition} of a matchable graph $G$ is:
\begin{enumerate}
\item a collection of Deming-BP subgraphs $\{ B_i\}_{i=1}^r$ of $G$ ($r\geq 0$);
\item a collection of Deming-$K_4$ subgraphs $\{K_i\}_{i=1}^{\ell}$ of $G$ $(\ell\geq 0)$;
\item a KE graph $R$ (possibly null); and such that
\item the vertices of these subgraphs partition $V(G)$.
\end{enumerate}
Figure~\ref{Deming-decomp-examp} depicts a Deming decomposition in
which each `part' has a single constituent.

It is worth emphasizing that any graph with a perfect matching has a 
Deming decomposition $\{\{ B_i\}_{i=1}^r, \{ K_j\}_{j=1}^{\ell},R\}$,
but a Deming decomposition is not necessarily unique. The input graph
$G$ is KE if and only if 
$r=\ell=0$. If $G$ is Egerv\'{a}ry, then it has no Deming-BP subgraphs
and $r=0$. 
A Deming decomposition is then a tool for investigating whether a
graph is Egerv\'{a}ry. We shall see that it yields necessary conditions
for whether a graph is Egerv\'{a}ry. We shall also see that this
decomposition can be extended to general graphs.  

\begin{center}
\begin{figure}
\includegraphics[width=9cm]{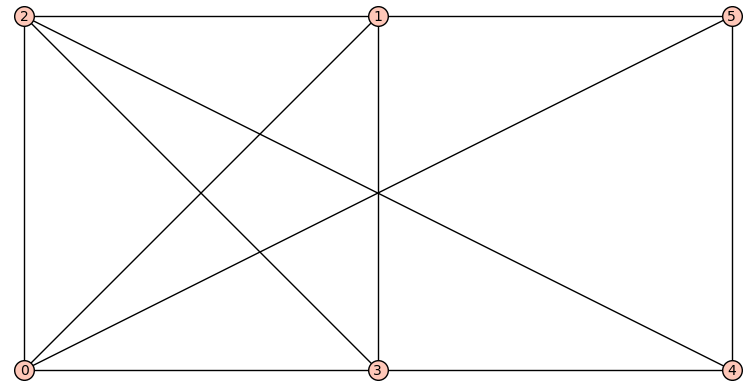}
\caption{An initial matching $\{\text{1--5, 3--4, 0--2}\}$ may differ
  from a Deming decomposition induced perfect matching 
$\{\text{0--2, 1--3, 4--5}\}$.} 
\label{k4notminimal}
\end{figure}
\end{center}

In computing a Deming decomposition, start by computing an initial
perfect matching $M$ (which won't necessarily
induce a perfect matching in each of the Deming subgraphs).
Consider, for instance, the graph in Figure~\ref{k4notminimal} and the perfect
matching $\{\text{1--5, 3--4, 0--2}\}$. A Deming decomposition consists of
the Deming-$K_4$ graph induced on the vertices $\{0,1,2,3\}$ and the
KE graph induced on the vertices $\{4,5\}$. So $M$ doesn't induce a
perfect matching on either subgraph. Of course, since each subgraph in a
Deming decomposition has a perfect matching, these define a perfect
matching $M'$ of the parent graph. We call such an $M'$ a
\textit{Deming decomposition induced perfect matching}.  
(See the second paragraph following the statement of 
Corollary~\ref{alpha_same_in_R_and_R'} for related remarks.)

For future reference, we collect several of our results on Deming
decompositions into a single theorem.

\begin{thm}\label{thm:dd}
If $G$ is a matchable graph with a Deming decomposition
$\{\{ B_i\}_{i=1}^r, \{ K_j\}_{j=1}^{\ell},R\}$, then:
\begin{enumerate}
\item $\alpha(B_i) = \nu(B_i)-1$, for $i \in [r]$;
\item $\alpha(K_j) = \nu(K_j)-1$, for $j \in [\ell]$;
\item $\alpha(R)=\nu(R)$;
\item\label{nu-adds} $\nu(G)=\sum_{i=1}^r \nu(B_i)+\sum_{j=1}^{\ell} \nu(K_j)+\nu(R)$; and
\item\label{alpha-sub-adds} $\alpha(G)\leq \nu(G)-(r+\ell)$.
\end{enumerate}
\end{thm}

\begin{center}
\begin{figure}
\includegraphics[width=7cm]{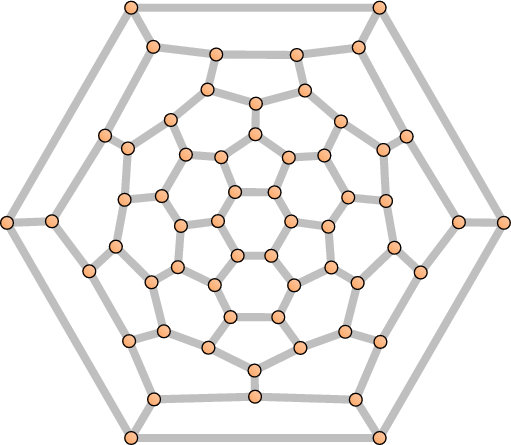}
\caption{The Buckminster Fullerene $C_{60}$ has 60 vertices. Of its $32$ 
  faces, $12$ are pentagons which are pairwise vertex disjoint.  One Deming decomposition
  consists of $6$ pairs of pentagons, each pair  joined by a single
  edge. Interestingly, $\alpha(C_{60})=24$, and the independence
  number of each of these $6$ Deming subgraphs is $4$.} 
\label{c60}
\end{figure}
\end{center}

Before moving to the next section, where we employ our basic results 
in further understanding Egerv\'{a}ry graphs, let us consider another
concrete example: the so-called `Buckminster Fullerene'; see
Figure~\ref{c60}.  While Theorem~\ref{thm:dd} (part~(\ref{nu-adds}))
shows that the matching number $\nu$ is additive across the parts of a 
Deming decomposition, part~(\ref{alpha-sub-adds}) guarantees only that the
independence number $\alpha$ is subadditive. Nevertheless, at least
for the graph $C_{60}$, the invariant $\alpha$ does behave additively
for the decomposition mentioned in Figure~\ref{c60}'s caption.

\section{Characterizing \& Recognizing Egerv\'{a}ry Graphs}
\label{sec:egervary}

Deming's Algorithm can be used to identify whether an input graph $G$
with a perfect matching $M$ is KE. If it is, the algorithm also produces a
maximum independent set $I$, which yields a certificate (in this case
$|M|+|I|$ must equal the order). In the case where $G$ is not KE, the
algorithm either produces a nice even $T$-subdivision, and this
provides a certificate that $G$ is not Egerv\'{a}ry. Deming's
Algorithm may also terminate with the production of a nice even
$K_4$-subdivision $F$.  While the production of $F$ provides a
certificate that $G$ is not KE, it may or may not be the case that $G$
is Egerv\'{a}ry. It remains an open question to find a certificate
that a graph $G$ is Egerv\'{a}ry.  
We record some conditions that may be relevant.

\begin{thm}\label{main}
A graph $G$ with a perfect matching is Egerv\'{a}ry if and only if $G$
is not an even $T$-subdivision and, for every edge $e$ not in all
perfect matchings of $G$, the graph $G-e$ is Egerv\'{a}ry. 
\end{thm}

\begin{proof}
Suppose first that $G$ is Egerv\'{a}ry. So $G$ does not have a nice
even $T$-subdivision. In particular, $G$ is not itself a
$T$-subdivision. If every edge is in every perfect matching of $G$,
then $G$ must consist of isolated edges and the condition that---for
every edge $e$ not in all perfect matchings of $G$, the graph $G-e$ is
Egerv\'{a}ry---is vacuously satisfied. Assume then that there is an
edge $e$ that is not in every perfect matching of $G$. Note that $G-e$
must then have a perfect matching and that $G-e$ is Egerv\'{a}ry: if
$G-e$ had a nice even $T$-subdivision, then $G$ would as well.  

Conversely, suppose that $G$ is not an even $T$-subdivision, and, for
every edge 
$e$ not in all perfect matchings of $G$, the graph $G-e$ is Egerv\'{a}ry. If
$G$ has a nice even $T$-subdivision $H$, then, since $G$ does not equal
$H$, the parent $G$ must contain an edge $e$ incident to some vertex
of $H$. Since $H$ is nice, it
has a perfect matching $M$ that can be extended to all of $G$. Since
$e\not\in M$, the edge $e$ is not in every perfect matching. But then
$G-e$ contains the nice even $T$-subdivision $H$, contradicting the
fact that $G-e$ is Egerv\'{a}ry. Therefore, $G$ contains no such $H$
and by Corollary~\ref{Egervary_char-even-T} is Egerv\'{a}ry.
\end{proof}

Given a matchable graph $G$ with a Deming decomposition 
$\{\{B_i\}_{i=1}^r, \{ K_j\}_{j=1}^{\ell},R\}$,  
one obvious necessary condition for $G$ to be Egerv\'{a}ry is that
there cannot be any Deming-BP subgraphs: if $G$ has a Deming-BP
subgraph, then it has a nice even $T$-subdivision and is not
Egerv\'{a}ry (by Corollary~\ref{Egervary_char-even-T}).  
We will establish other necessary conditions. 

One other easy-to-see necessary condition is that each Deming subgraph
itself is Egerv\'{a}ry. Yet another is that there can be no edge
between vertices of different Deming subgraphs: these have spanning
even $K_4$ subdivisions and if there is such an edge, then a cycle
from each subgraph can be paired to construct a nice even
$T$-subdivision.  

The following theorem resembles Theorem~\ref{main},
specified to the Deming-$K_4$ subgraph case. In fact, Deming-$K_4$
subgraphs have 
significant structure and this specification utilizes that
structure. In particular, Deming-$K_4$ subgraphs $D$, by definition,
do not have 
nice even $K_4$-subdivisions with fewer vertices: if one did, there
would be an edge $xy$ of $D$ not spanned by that (smaller) subdivision
and 
then, by definition $D-\{x,y\}$ would both be KE and have a nice
even $K_4$-subdivision, a contradiction. In this sense, Deming-$K_4$ 
subgraphs are 
minimal non-KE graphs with nice even $K_4$-subdivisions. The
proof of this next result depends on this fact essentially.

\begin{thm}
A Deming-$K_4$ graph $K$ with a spanning even $K_4$-subdivision $F$ is 
Egerv\'{a}ry if and only if for every edge $e$ of $F$, either: (1)
$K-e$ is KE, or (2) $K-e$ itself has a spanning even $K_4$-subdivision
and is Egerv\'{a}ry. 
\end{thm}

\begin{proof}
First assume that $K$ is Egerv\'{a}ry. 
Given an edge $e$ of $F$, we'll show that 
either (1) or (2) holds.
Suppose that $K-e$ is not KE and thus contains either a nice
even $T$-subdivision or a nice even $K_4$-subdivision. First consider
the case where $K-e$ contains a nice even $T$-subdivision $H$. But
then $H$ is  a nice even $T$-subdivision within $K$, 
contradicting the
assumption that $K$ is Egerv\'{a}ry. So it must be that $K-e$ contains
a nice even $K_4$-subdivision $F'$. If $F'$ does not span $K-e$,
then---since $F'$ is a nice subgraph of $K-e$---the vertices of 
$K-e$ that are not in $F'$ are perfectly matched by 
some edges of $K-e$. 
If $xy$ is a perfect matching edge between two vertices 
not in $F'$, then $xy$ must be a perfect
matching edge of $K$. Since $K$ is a Deming-$K_4$ subgraph, it follows
that $K-\{x,y\}$ is KE. But $K-\{x,y\}$ contains a nice even
$K_4$-subdivision $F'$, contradicting the fact that $K-\{x,y\}$ is
KE. So it must be that $F'$ spans $K-e$. Since $K$ is Egerv\'{a}ry,
$K-e$ must also be Egerv\'{a}ry, and we have that  $K-e$ itself has a
spanning even $K_4$-subdivision and is Egerv\'{a}ry.  

For the converse, suppose that $K$ is not Egerv\'{a}ry.
By Corollary~\ref{Egervary_char-even-T}, then $K$ has a nice even
$T$-subdivision $H$. 
The subgraph $F$ must contain an edge $e$ not in $H$, and 
we can assume that $e$ is incident to a vertex of $H$. 
The graph $K-e$ cannot be KE or Egerv\'{a}ry as $H$ is a nice even
$T$-subdivision of $K-e$. But then neither (1) nor (2) holds. 
\end{proof}

While we are lacking sufficient conditions for a graph with a perfect
matching to be Egerv\'{a}ry, we do have a few more necessary
conditions, that, accumulated, will suggest a characterization
conjecture (see Conjecture~\ref{conj:Egervary-char}).

\begin{thm}
If $G$ is an Egerv\'{a}ry graph with a perfect matching and
corresponding Deming decomposition $\{\{ K_j\}_{j=1}^{\ell},R\}$, and
$K$ and $K'$ are Deming-$K_4$ subgraphs, 
then no edge is incident to both a vertex in $K$ and one in $K'$. 
\end{thm}

\begin{proof}
Denote a perfect matching in $G$ by $M$, which we may assume induces
perfect matchings $M_1$ in $K$  and $M_2$ in $K'$. We will identify $K$ and $K'$ 
with their spanning $K_4$-subdivisions---the `extra' edges (edges not
in these spanning subgraphs) will play no role here in our
considerations. 

Let $v_1,v_2,v_3,v_4$ be the corner
vertices of $K$ and $w_1,w_2,w_3,w_4$ be the corner vertices of
$K'$. Let $P_{i,j}$ be the odd paths connecting corner vertices in
$K$ and $Q_{i,j}$ be the corresponding paths in $K'$. Suppose that
$x\in V(K)$ and $y\in V(K')$, with $x$ adjacent to $y$. We will show
that then $G$ contains a nice even subdivision of $T$---contradicting the
fact that $G$ is Egerv\'{a}ry. 

In the first case, assume that both $x$ and $y$ are corner
vertices. We may assume that $x=v_1$, $y=w_1$, that  $P_{1,2}$ and
$P_{3,4}$ are $M_1$-saturated paths,  and that $Q_{1,2}$ and $Q_{3,4}$
are $M_2$-saturated  paths. Thus, the
cycle formed from $P_{2,3}$ followed by $P_{3,4}$ and then $P_{4,2}$ is a
blossom with tip $v_2$. Similarly, the cycle formed from
$Q_{2,3}$ followed by $Q_{3,4}$ and then $Q_{4,2}$ is a blossom with
tip $w_2$. The $M$-saturated $(v_2,w_2)$-path 
formed from $P_{2,1}$ followed by the edge $v_1w_1$ (which is $xy$) and then
$Q_{1,2}$ is an odd path between blossom tips. Thus, the two blossoms and
this path form an even subdivision of $T$. Now edges in $M$ are
incident to either zero or two vertices in this subgraph, and since $M$ is
perfect, we have a nice subgraph. 

In the second case, assume that $x$ is \emph{not} a corner vertex and
that $y$ \emph{is} a corner vertex. We may assume that $x$ is
interior to $P_{1,2}$ and that $y=w_1$; we continue to assume 
that $P_{1,2}$ and $P_{3,4}$ are $M_1$-saturated  and that $Q_{1,2}$ and $Q_{3,4}$
are $M_2$-saturated. The key point here is that $x$ is an odd distance
from exactly one of $v_1$, $v_2$. Assume that $x$ is an odd distance
from $v_2$. Again the cycle formed from $P_{2,3}$ followed by
$P_{3,4}$ and then $P_{4,2}$ is a blossom with tip $v_2$, and the
cycle formed from $Q_{2,3}$ followed by $Q_{3,4}$ and then $Q_{4,2}$
is a blossom with tip $w_2$. The path from $v_2$ to $x$
to $w_2$ formed from the $(v_2,x)$-segment of $P_{2,1}$, followed by
the edge $xy$ and then $Q_{1,2}$ is an odd $M$-saturated path
joining blossom tips, and so we again have a nice even subdivision of $T$ as
a subgraph of $G$. 

Finally, consider the case where neither $x$ nor $y$ is a corner
vertex. Here, we may assume that $x$ is interior to $P_{1,2}$ and that
$y$ is interior to $Q_{1,2}$. We may further assume that $P_{1,3}$
and $P_{2,4}$ are $M_1$-saturated paths and that $Q_{1,3}$
and $Q_{2,4}$ are $M_2$-saturated paths. Similarly to the preceding
case, we may here assume that $x$ is a (positive) even distance from
$v_1$ and that $y$ is a (positive) even distance from $w_1$. Now the
cycle formed from $P_{3,4}$ followed by $P_{4,2}$ and then $P_{2,3}$
is a blossom $\EuScript{B}_1$ with tip $v_3$, and the
cycle formed from $Q_{3,4}$ followed by $Q_{4,2}$ and then $Q_{2,3}$
is a blossom $\EuScript{B}_2$ with tip $w_3$. The odd $M$-saturated
path $\EuScript{P}$ joining the tips $v_3$, $w_3$ here takes the form
$P_{3,1}$ followed by the $(v_1,x)$-segment of $P_{1,2}$, then by the
edge $xy$, then by the $(y,w_1)$-segment of $Q_{1,2}$, and finally by 
$Q_{1,3}$. Once again, we've revealed a nice even subdivision
$\EuScript{B}_1\cup\EuScript{P}\cup\EuScript{B}_2$ of $T$ in $G$, the
desired contradiction.
\end{proof}

\begin{thm}
\label{thm:no_alt_vv-path}
If $K$ is a Deming-$K_4$ graph with perfect matching $M_K$, $R$ is a
vertex-disjoint KE graph with perfect matching $M_R$,  and $G$ is an
Egerv\'{a}ry graph formed by adding some edges between $K$ and $R$,
then no vertex $v\in V(K)$ admits an $M_R$-alternating path to itself. 
\end{thm}

\noindent
\textbf{Remark:} Technically, we should use `cycle' instead of `path'
in the assertion because the walk under consideration would be
closed. But, as seen in the proof, it would also have odd length,
hence vacuously not be $M_R$-alternating. The point of the result is
not to rule out such a cycle for this (trivial) parity reason, rather
to rule it out for the structural reason explored in the proof.

\begin{proof}
Let $v$ be a vertex in $K$ and $F$ be $K$'s spanning even
$K_4$-subdivision containing $M_K$. We need to show that there is no
$M_R$-alternating path from $v$ to itself. For a contradiction,
suppose that there is such a path, say $e_1,e_2,...,e_r$.  As
$M_R$-edges lie in $R$ and $v\not\in V(R)$, the edges $e_1$ and
$e_r$ are not in $M_R$. So this path (really cycle) has odd length
and hence defines a blossom $\EuScript{B}_1$ in $G$. Since $M_K$ is
a perfect matching of $K$, the vertex $v$ must be incident to an
$M_K$-edge $vw$ on one of the paths $P_{i,j}$ defining $F$---we can
assume that it's $P_{1,2}$ and that $w$ is closer to $v_2$ than to
$v_1$. There are two cases to consider. In the first case, the
distance from $w$ to $v_2$ along $P_{1,2}$ is even. Here, the last
edge of the $(w,v_2)$-segment of $P_{1,2}$ is necessarily a matching
edge. So the cycle formed by concatenating $P_{2,3}$, $P_{3,4}$, and
$P_{4,2}$ is a blossom $\EuScript{B}_2$, and the blossom
$\EuScript{B}_1$, followed by the $(v,v_2)$-segment of $P_{1,2}$ (an
odd path), and followed by $\EuScript{B}_2$ is a nice even
$T$-subdivision. But then $G$ is not Egerv\'ary, which is a
contradiction. In the second case, the distance from $w$ to $v_2$
along $P_{1,2}$ is odd, so the last edge on the $(w,v_2)$-segment of
$P_{1,2}$ is not in $M_K$. But $v_2$ is saturated by $M_K$; so
assume that the first (and hence the last) edge of the path
$P_{2,3}$ is in $M_K$. Now the concatenation of $P_{3,4}$,
$P_{4,1}$, and $P_{1,3}$ is a blossom $\EuScript{B}_2$ with blossom
tip $v_3$.  So the blossom $\EuScript{B}_1$, followed by the path
from $v$ to $v_2$ along $P_{1,2}$, followed by the path $P_{2,3}$
(an odd path), and finally followed by $\EuScript{B}_2$ is a nice
even $T$-subdivision. But then again we reach the contradiction that
$G$ is not Egerv\'ary.
\end{proof}

\begin{thm}\label{thm:kextra}
If $K$ is an  Egerv\'{a}ry Deming-$K_4$ graph with a spanning even
$K_4$-subdivision $F$ and $vw\in E(K)\smallsetminus E(F)$, then $vw$ lies
in a perfect matching of $K$. 
\end{thm}

\begin{proof}
Call the edges of $K$ not in $F$ \textit{extra} edges. 
By Corollary~\ref{cor:all_K4edges_in_PM}, the conclusion holds
when $K$ has $0$ extra edges, for such $K$ are just even
$K_4$-subdivisions. Assume, then, that for a fixed integer $k>0$, the
conclusion holds for Egerv\'{a}ry 
Deming-$K_4$ graphs with fewer than $k$ extra edges, and suppose that
$K$ has $k$ extra edges $e_1,e_2,\ldots,e_k$. 

We need to show that each of these edges is in a perfect matching of $K$,
and when $k>1$, this follows immediately by induction. For if 
$1\leq i\leq k$, then the subgraph $K'=K-e_i$ is Egerv\'{a}ry (because 
$K$ is Egerv\'{a}ry), has $F$ as a spanning even
$K_4$-subdivision, and has $k-1>0$ extra edges. So all these extra edges
are in perfect matchings of $K'$ (hence of $K$), and since $i$ was
arbitrary, the conclusion follows. 

Indeed, the crux of this proof is when $k=1$, i.e., when $e=e_1$ is
the only extra edge (so $K$ consists of an even $K_4$-subdivision
with one extra edge $e=vw$).  

There are three main cases to consider: (1) when $e$ is
incident to a pair of vertices on a single $P_{i,j}$ path of $F$ (we can
assume $P_{1,2}$); (2) when $e$ is incident to vertices on two
`incident' paths of $F$ (we'll assume $P_{1,2}$ and $P_{1,3}$); and (3)
when $e$ is incident to vertices on `opposite' paths of $F$ (we'll
assume $P_{1,2}$ and $P_{3,4}$). In each case, we'll show either that
the case is impossible---that the existence of such an edge violates
the assumption that $K$ is Egerv\'{a}ry---or that $K$ contains a
perfect matching containing $e$. 

In case (1), we assume that $e$'s endpoints $v$, $w$ on $P_{1,2}$
have $v$ closer to $v_1$. We first argue that the number of edges
from $v$ to $w$ on $P_{1,2}$ must be odd.  If this number were even, we
could find a nice even $T$-subdivision in $K$, contrary to 
$K$ being Egerv\'{a}ry. First note that exactly one of
the subpaths from $v$ to $v_1$ or $w$ to $v_2$ of $P_{1,2}$ is
even. Let us assume that the first of these (from $v$ to $v_1$) is
even (the other possibility is similar). Since the number of edges
from $v$ to $w$ on $P_{1,2}$ is even, the cycle formed by these edges
together with the edge $vw$ is odd. Next, there is an odd path from
$v$ to $v_1$ (along $P_{1,2}$) to $v_4$ (using all of
$P_{1,4}$). Finally, the paths $P_{4,2}$, followed by $P_{2,3}$, and
then $P_{3,4}$ form an odd cycle. We have now formed our desired even
$T$-subdivision. And it's nice because the omitted vertices of $K$ are
those strictly between $w$ and $v_2$ along $P_{1,2}$ and those
strictly between $v_1$ and $v_3$ along $P_{1,3}$, all of which can be
perfectly matched along these paths.

Now we have that the number of edges from $v$ to $w$ along $P_{1,2}$
is odd.  Thus, the $(v,v_1)$- and $(w,v_2)$-segments of $P_{1,2}$ have
the same parity. We consider the subcase when these segments are both
of even length; the odd case is similar.  Here, we find a perfect
matching of $K$ consisting of $vw$, independent edges between $v$ and
$w$ on $P_{1,2}$, independent edges between $v$ and $v_1$---including
one that is necessarily incident to $v_1$---independent edges between
$w$ and $v_2$---including one that is necessarily incident to
$v_2$---independent edges on $P_{3,4}$ including one incident to $v_3$
and one incident to $v_4$ that saturate all vertices on this path, and
lastly independent edges saturating all subdivision vertices on
$P_{1,4}$, $P_{2,3}$, $P_{1,3}$, and $P_{2,4}$.

In case (2), we assume that $e$'s endpoint $v$ is on
$P_{1,2}$ and $w$ is on  $P_{1,3}$. We shall argue momentarily that
the number $N$ of 
edges on the path $P$ from $v$ to $v_1$ along $P_{1,2}$ followed by
$v_1$ to $w$ on $P_{1,3}$ can be assumed odd. In this case, we
construct a perfect 
matching of $K$ starting with $vw$ together with independent edges
between $v$ and $w$ on $P$. 
Either the number of edges from $v$ to $v_2$ on $P_{1,2}$
is even and the number of edges from $w$ to $v_3$ on $P_{1,3}$ is odd,
or vice versa.
First assume the former.
Then it is possible to extend the matching to saturate
all remaining vertices on $P_{1,2}$ including $v_2$ and every vertex
of $P_{1,3}$ except $v_3$. The matching can be further extended with
independent edges along $P_{3,4}$ that saturate all these vertices and
finally covering the subdivision vertices of $P_{2,4}$, $P_{1,4}$, and
$P_{2,3}$. This yields a perfect matching of $K$ which includes the edge
$vw$.  The latter subcase---when the number of edges from $v$ to $v_2$
on $P_{1,2}$ is odd and the number of edges from $w$ to $v_3$ on
$P_{1,3}$ is even---can be addressed similarly.

We turn to the case when $N$ (of the preceding paragraph) is even.
There are two subcases
depending on whether the number $N'$ of edges from $v$ to $v_2$ along
$P_{1,2}$ is odd or even (note that the number of edges from $w$ to
$v_3$ along $P_{1,3}$ must have the same parity as $N'$). Consider the
subcase 
when $N'$ is odd. Since $N$ is even,  the path $P$
together with the edge $vw$ forms an odd cycle $C$ containing $v$. 
From there, the $(v,v_2)$-segment of $P_{1,2}$ is also odd (because
$N'$ is odd), and this path attaches at $v_2$ to the odd cycle formed
by the paths $P_{2,4}$, $P_{4,3}$, and 
$P_{3,2}$. Thus we observe an even $T$-subdivision in $K$.
The remaining vertices (not
belonging to $T$) are the ones strictly between $w$ and $v_3$ on
$P_{1,3}$ (an even number) and the ones strictly between $v_1$ and
$v_4$ on $P_{1,4}$ (also an even number). Since these can all
be saturated by a matching, we see that $T$ is nice, which is impossible.

Within case (2), it remains to consider the subcase when $N'$ is even,
and here, we construct another perfect matching $M$ of $K$ containing the
edge $vw$. With $N$ still being even, we see that the odd
cycle $C$ from above persists. Thus, we can initialize $M$ with edges
including $vw$ and saturating every vertex of $C$ except $v_1$. With 
$N'$ being even, the $(v,v_2)$-segment of $P_{1,2}$ and the 
$(w,v_3)$-segment of $P_{1,3}$ are both even. Hence we can extend $M$
with independent edges along these segments so that $v_2$ and $v_3$
are both saturated, along with the internal vertices up to but not
including $v$ and $w$ (as these were already covered).
To reach a perfect matching, it remains to include
independent edges along $P_{1,4}$ that saturate all these vertices and
finally cover the subdivision vertices of $P_{2,4}$, $P_{3,4}$, and
$P_{2,3}$; this is possible because these paths are all odd.

The last case (3) is in a sense the easy case because none of the
subcases lead to contradictions. Here, the edge $e=vw$ is incident to
vertices on `opposite' paths (we will 
assume $P_{1,2}$ and $P_{3,4}$, with $v$ on the first and $w$ on the
second). In this case, $e$ can always be included in a perfect
matching of $K$. The vertex $v$ must be either at an odd distance from
$v_1$ on $P_{1,2}$ and an even distance from $v_2$ or at an even
distance from $v_1$ on $P_{1,2}$ and an odd distance from $v_2$. We
assume the former and omit the latter (which can be addressed
similarly). 

In the first subcase, the distance from $w$ to $v_4$ on $P_{3,4}$ is
odd. Then the edge $e$ is included in the even cycle consisting of the
path from $v$ to $v_1$ on $P_{1,2}$, followed by $P_{1,4}$, followed
by the path from $v_4$ to $w$ on $P_{3,4}$, and concluding with the
edge $e$ (connecting vertex $w$ back to starting point $v$). The
remaining vertices consist of 
three paths, each with an even number of
vertices. The observed even cycle can be perfectly matched with a
matching containing $e$, and the three remaining paths each can be
perfectly matched, all together yielding a perfect matching of $K$
containing $e$.  

In the remaining subcase, the distance from $w$ to $v_4$ on $P_{3,4}$
is even. Here, the edge $e$ is included in the even cycle
consisting of the path from $v$ to $v_1$ on $P_{1,2}$, followed by
$P_{1,3}$, followed by the path from $v_3$ to $w$ on $P_{3,4}$,
and finishing with the edge $e$ (connecting $w$ back to the starting
point $v$). The remaining vertices again comprise
three odd-length paths
(i.e., each with an even number of vertices). Again, the described
even cycle can be perfectly 
matched with a matching containing $e$, and the remaining paths can
each be perfectly matched, yielding a perfect matching of $K$
containing $e$. 
\end{proof}

\begin{cor}\label{cor:kperfect}
If $K$ is a Deming-$K_4$ graph that is Egerv\'{a}ry, then every edge
of $K$ is in some perfect matching of $K$. 
\end{cor}

The following theorem simply summarizes our previous results. 
The proof of condition (\ref{no-alt-path-item}) therein is analogous
to that of (the more difficult) Theorem~\ref{thm:no_alt_vv-path}.
The graph $G[V(K)\cup V(R)]$ is the subgraph of $G$ induced on the
union of the vertices of subgraphs $K$ and $R$, that is, the graphs
$K$ and $R$ together with any edges between them in $G$. 

\begin{thm}\label{thm:egervary}
Let $G$ be a graph with a perfect matching $M$ and a Deming
decomposition $\{\{ B_i\}_{i=1}^r, \{ K_j\}_{j=1}^{\ell},R\}$  (with
respect to $M$).  If $G$ is Egerv\'{a}ry then: 
\begin{enumerate}
\item there are no Deming-BP subgraphs (that is, $r=0$);
\item each Deming-$K_4$ subgraph $K_j$ is Egerv\'{a}ry;
\item every edge in each Deming-$K_4$ subgraph $K_j$ is in a perfect matching of $G$;
\item the Deming-$K_4$ subgraphs $K_j$ are independent (that is, there are no edges in $G$ incident to more than one Deming-$K_4$ subgraph);
\item for each Deming-$K_4$ subgraph $K_j$, the graph $G[V(K_j)\cup V(R)]$ is Egerv\'{a}ry; and 
\item\label{no-alt-path-item} for the matching $M_R$ consisting of $M$
  restricted to $R$, 
  there are no $M_R$-alternating paths between any pair of vertices
  $v,w$ from different Deming-$K_4$ subgraphs. 
\end{enumerate}
\end{thm}

We conjecture that this necessary condition is also sufficient---that
is, this is a characterization of Egerv\'{a}ry graphs---and also that
the independence number of $G$ is the sum of the independence numbers
of these Deming decomposition subgraphs. We record these formally.

\begin{conj}
\label{conj:Egervary-char}
Let $G$ be a graph with a perfect matching and a Deming decomposition 
$\{\{ B_i\}_{i=1}^r, \{ K_j\}_{j=1}^{\ell},R\}$.  The $G$ is
Egerv\'{a}ry if and only if:
\begin{enumerate}
\item there are no Deming-BP subgraphs (that is, $r=0$);
\item each Deming-$K_4$ subgraph $K_j$ is Egerv\'{a}ry; 
\item every edge in each Deming-$K_4$ subgraph $K_j$ is in a perfect matching of $G$;
\item the Deming-$K_4$ subgraphs $K_j$ are independent;
\item for each Deming-$K_4$ subgraph $K_j$, the graph $K_j+R$ is Egerv\'{a}ry; and 
\item there are no $M_R$-alternating paths between any pair of
  vertices $v$, $w$ from different Deming-$K_4$ subgraphs. 
\end{enumerate}
\end{conj}

\begin{conj}
If $G$ is Egerv\'{a}ry with Deming decomposition 
$\{\{K_j\}_{j=1}^{\ell},R\}$, then 
$\alpha(G)=\sum_{j=1}^{\ell}\alpha(K_j)+\alpha(R)$ (and thus $\alpha(G)=\nu(G)-\ell$). 
\end{conj}

\section{Constructing Egerv\'{a}ry Graphs}
\label{sec:constructing}

Every KE graph is Egerv\'{a}ry
(Corollary~\ref{cor:kayll-KE-implies-eger}). But there are also non-KE
graphs which 
are Egerv\'{a}ry. We give examples of some constructions. These
examples are interesting not only because they provide infinite families
of non-KE Egerv\'{a}ry graphs---but also because they are not
obviously derived from Deming subgraphs with the addition of extra
edges.

The graph $K_4$ is the smallest instance of what we call a
\textit{weak wheel}, i.e., a
cycle $C$ with three or more nodes, a node $w$ not in $C$, and
a path (\textit{spoke}) from each node of $C$ to $w$ such that these
spokes are disjoint except for $w$. 

\begin{thm}
\label{weak-wheels-Eger}
Every weak wheel with a perfect matching is Egerv\'{a}ry.
\end{thm}

\begin{proof}
A weak wheel cannot have a disjoint pair of odd cycles.
\end{proof}

\begin{figure}
\begin{tabular}{ccc}
\begin{subfigure}{0.35\textwidth}
\includegraphics[width=0.9\linewidth]{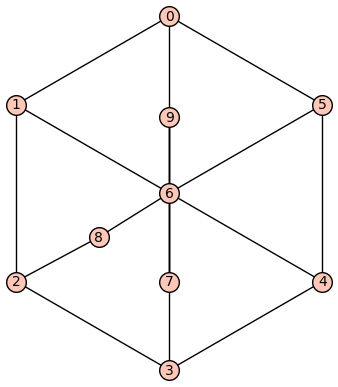}
\caption{$G_1$}
\end{subfigure}
&
\begin{subfigure}{0.35\textwidth}
\includegraphics[width=0.9\linewidth]{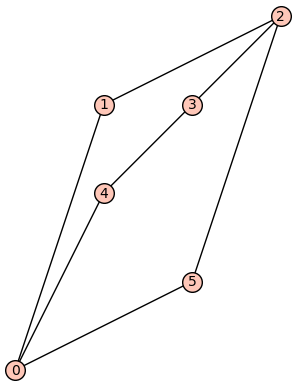}
\caption{$G_2$}
\end{subfigure}
&
\begin{subfigure}{0.35\textwidth}
\includegraphics[width=0.9\linewidth]{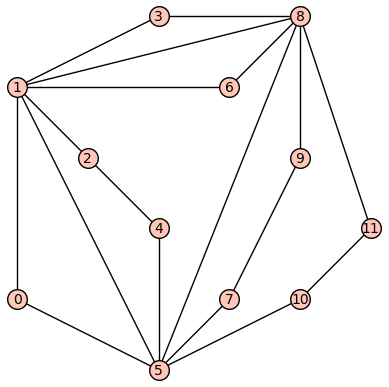}
\caption{$G_3$}
\end{subfigure}
\end{tabular}
\caption{Graph $G_1$ is a \textit{weak wheel}, $G_2$ is a \textit{weak
    banana}, and $G_3$ is a \textit{bracelet}. They all have perfect
  matchings, so by
  Theorems~\ref{weak-wheels-Eger}--\ref{bracelets-Eger}, resp., are
  Egerv\'{a}ry.}
\label{fig:Jacks-examples}
\end{figure}

A graph is a \textit{weak banana} when it consists of
two nodes, $v$, $w$, together with three or more internally disjoint
$vw$-paths. Notice that a weak banana fails to be bipartite
exactly when at least two of its $vw$-paths are of opposite
parity.

\begin{thm}
\label{weak-bananas-Eger}
Every weak banana with a perfect matching is Egerv\'{a}ry.
\end{thm}

\begin{proof}
A weak banana cannot have a disjoint pair of odd cycles.
\end{proof}

A \textit{bracelet} consists of three distinct nodes $u$, $v$, and
$w$, together with three weak bananas, respectively joining $u$
to $v$, $v$ to $w$, and $w$ to $u$. 

\begin{thm}
\label{bracelets-Eger}
Every bracelet with a perfect matching is Egerv\'{a}ry.
\end{thm}

\begin{proof}
A bracelet cannot have a disjoint pair of odd cycles.
\end{proof}

Figure~\ref{fig:Jacks-examples} depicts one example each of a weak wheel,
a weak banana, and a bracelet.

\begin{figure}
\includegraphics[width=350pt]{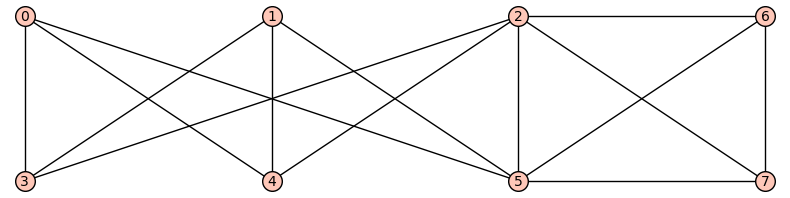}
\caption{a bipartite extension $K_{4} : K_{3,3}$}
\label{fig:bip-extension}
\end{figure}

A \textit{bipartite extension} of a graph $G$ is obtained by 
replacing an edge $e=vw$ of $G$ by a connected bipartite graph
$H$ so that $v$ is identified with a node in one part of $H$ and
$w$ with a node in the other; we use
$G:H$ to denote the resulting graph. 
See Figure~\ref{fig:bip-extension} for the bipartite extension $K_4 : K_{3,3}$,
in which the edge 2--5  of $K_4$ has been deleted, the vertex $2$ is
the identification of the vertex $2$ from $K_4$ and one from $K_{3,3}$, and
the vertex $5$ is the identification of the vertex $5$ from $K_4$ and one
from the opposite part of $K_{3,3}$. It's an exercise to confirm that
a bipartite extension $G:H$ of a matchable graph $G$ is itself
matchable only if the bipartite graph $H$ has parts of equal size.

The following assertion is immediate.
\begin{thm}
\label{bip-exts-preserve-Eger}
If a matchable graph $G$ contains no disjoint pair
of odd cycles, then every matchable bipartite extension of $G$ is Egerv\'ary.
\end{thm}

\begin{proof}
Under the hypotheses, a bipartite extension of $G$ cannot contain a
disjoint pair of odd cycles. 
\end{proof}

\noindent
Theorem~\ref{bip-exts-preserve-Eger} provides a way to
generate several infinite families of Egerv\'{a}ry graphs that
are not K-E: every bipartite extension of a weak wheel (or a
non-bipartite weak banana, or a bracelet) that has a perfect
matching is non-KE  but Egerv\'{a}ry. 

Notice that the Egerv\'{a}ry graphs of
Theorems~\ref{weak-wheels-Eger}--\ref{bip-exts-preserve-Eger} 
all can be obtained by
starting with a graph $G$ that contains no disjoint pair of
(either even or odd) cycles, then replacing some edges of $G$
by paths to get some odd cycles, and finally taking some
bipartite extensions.

Here is another way to produce new Egerv\'{a}ry graphs from old.

\begin{thm}
\label{Egervary-from-Eger-plus-KE}
If $H$ is Egerv\'{a}ry, $R$ is a matchable KE graph vertex-disjoint
from $H$, the set $I$ is a maximum independent set of $R$, and $G$ is formed by adding any edges from $H$ to $N=N(I)$, then $G$ is Egerv\'{a}ry. 
\end{thm}

\begin{proof}
We proceed by contradiction. 
If $G$ is not Egerv\'{a}ry, then $G$ contains a nice, even T-subdivision $S$.  
As $H$ is Egerv\'{a}ry, it doesn't contain $S$, and likewise the KE graph
$R$ doesn't contain $S$. So $V(S)$ meets both $V(H)$ and $V(R)$; let us
call edges of $S$ with one end in $V(H)$ and one end in $V(R)$ 
\DEF{crossing edges}.

Because $H$ is Egerv\'{a}ry (hence matchable), $R$ is matchable, and $S$ is
nice, $G$ contains a perfect matching $M$ with $M_S=M\cap E(S)$ a
perfect matching of $S$ and $M_R=M\cap E(R)$ a perfect matching of
$R$. Notice that every perfect matching of $G$ necessarily matches $I$
into $N$ because $I$ is independent and $G$ contains no edges from $I$
to $H$. In particular, $M$ matches $I$ into $N$. Also, $M$ matches $N$
into $I$ because $M_R$ is a perfect matching of $R$.

Now consider a crossing edge $xy$ of $S$ with $x\in V(R)$ and  
$y\in V(H)$. Then $x\in N$, and since $M$ matches $N$ into $I$, the edge
$xy$ is not in $M$. We distinguish two cases.

\noindent
Case~(i): $xy$ is in a cycle $C$ of $S$. 

As $M_S$ is a perfect matching of $S$ and $xy\not\in M$, the next edge
$xz$ of $C$, starting at $x$ and moving away from $y$, is in $M$, and,
knowing $M$ matches $N$ into $I$, we see that $z\in I$. Thus, the next
vertex $w$ of $C$ after $z$ has to be in $N$ and so cannot be $y$ (which
is in $V(H)$, not $N$). Continuing around $C$ from $w\in N$, we can
argue similarly as we did from $x\in N$ to see that the vertices of $C$
alternate between $N$ and $I$. But $C$ is an odd cycle, so when we
finally return to $y$, say, along an edge $uy$, the vertex $u$ must be
in $I$, and this forces $y$ to be in $N$. This is a contradiction
because $y\in V(H)$ and $V(H)\cap N=\varnothing$.

\noindent
Case~(ii): $xy$ belongs to the path $P$ joining the two cycles of $S$. 

As in Case~(i), the next vertex $z$ of $P$ moving away from $y$ must lie
in $I$ and the edge $xz$ must be in $M$. Continuing along $P$ from $z\in
I$, the edges must alternate out and in of $M$ because $M_S$ is a
perfect matching of $S$. And because $M$ perfectly matches $I$ with $N$,
the vertices $x,z\ldots$ of $P$ must lie alternately in $N$ and $I$. The
structure of the perfect matching $M_S$ within $S$ puts the blossom tip
$t$ reached along $P$ in this direction incident with an $M$-edge of
$P$. Thus, the $N$-$I$ alternation implies that $t$ belongs to
$I$. Being a tip, $t$ lies in one odd cycle $C'$ of $S$. With $t\in I$,
we see that one neighbor $v$ of $t$ on $C'$ must belong to $N$. Now
again, with $N$ matched to $I$ under $M$, the next vertex adjacent to
$v$ on $C'$ must lie in $I$. As in Case~(i), the vertices of $C'$ must
alternate between $I$ and $N$ as we traverse $C'$ starting a sequence $t\in
I$, $v\in N\ldots$. Thus, with $C'$ being odd, the other neighbor of $t$ on
$C'$ lies in $I$ together with $t$, which contradicts the independence
of $I$.
\end{proof}

\section{Independence Structure}
\label{sec:structure}

Some structure of maximum independent sets was given in \cite{Lars11}:
here we show that every graph can be decomposed into a KE graph and
a 2-bicritical graph with certain attractive properties. A Deming
decomposition of the (unique) 2-bicritical subgraph extends this
earlier decomposition. It is worth noting that
Pulleyblank~\cite{Pull79}  proved that
almost all graphs are 2-bicritical. 

Following a characterization in \cite{Pull79}, we call
a graph $G$ \textit{$2$-bicritical} if and only if $|N(I)|>|I|$ for
every nonempty independent set $I$ of $G$.  An
independent set $I$ is \textit{critical} if every independent set
$J$ satisfies the relation $|N(I)|-|I|\geq |N(J)|-|J|$. 
The neighbors $N(I)$ of a critical
independent set $I$ can necessarily be matched into $I$ (else it can
be shown that $I$ is not critical). A \textit{maximum critical
  independent set} is a critical independent set of largest
cardinality. The following theorem implies that a graph is
2-bicritical if and only if it has no nonempty critical independent
set. We can use this fact to provide a certificate that a graph is
2-bicritical. 

\begin{thm}[Larson, \cite{Lars11}]\label{idt}
For every graph $G$, there is a unique set $X\subseteq V(G)$ such that:
\begin{enumerate}
\item $\alpha(G)=\alpha(G[X])+\alpha(G[X^c])$;
\item $G[X]$ is a K\H{o}nig-Egerv\'ary graph;
\item $G[X^c]$ is 2-bicritical; and
\item each maximum critical independent set $J_c$ of $G$ satisfies $X=J_c\cup N(J_c)$.
\end{enumerate}
 \end{thm}
 
A Deming decomposition of $G[X^c]$ then extends this decomposition of
$G$.  Note that if $G$ has a perfect matching, it immediately follows
that the vertices of every maximum critical independent set $J_c$
must be matched to $N(J_c)$ and $|J_c|=|N(J_c)|$; it also necessarily
follows that $G[X^c]$ has a perfect matching. 
 
 \begin{prop}\label{prop:critical}
If $G$ has a critical independent set $I$ with $|I|=|N(I)|$, then $G$
is Egerv\'{a}ry if and only if $G-I-N(I)$ is Egerv\'{a}ry. 
\end{prop}

\begin{proof}
For the necessity, suppose that $G-I-N(I)$ is not Egerv\'{a}ry. By 
Corollary~\ref{Egervary_char-even-T}, then
this graph has a nice even subdivision $T_0$ of $T$. Because $I$ is
perfectly matched with $N(I)$, the subgraph $T_0$ is also nice in
$G$, and so $G$ is likewise not Egerv\'{a}ry.

Conversely, suppose that $G$ is not Egerv\'{a}ry. 
For a contradiction, suppose that $G-I-N(I)$ is Egerv\'{a}ry. 
With an eye toward applying Theorem~\ref{Egervary-from-Eger-plus-KE},
now let $H=G-I-N(I)$ and $R=G[I\cup N(I)]$ (the induced
subgraph). Then the triple $(H,R,I)$ satisfies the hypotheses of 
Theorem~\ref{Egervary-from-Eger-plus-KE}, and the present graph $G$ is
formed by adding some (maybe zero) edges from $H$ to $N(I)$ and so, by
the same theorem, is Egerv\'{a}ry. This is a contradiction.

Therefore, $G-I-N(I)$ is not Egerv\'{a}ry, and the proof is complete.
\end{proof}

The class of $\alpha$-critical graphs was mentioned earlier
(see the paragraph preceding Corollary~\ref{cor:alphacriticalk4}); in
particular, a  Deming-$K_4$ graph is $\alpha$-critical if and only if
it is an even  subdivision of $K_4$
(Corollary~\ref{cor:alphacriticalk4}). It is not  the case that
$\alpha$-critical graphs admit nontrivial 
decompositions with independence additivity.  
 
 \begin{prop}\label{prop:V1V2}
 If $G$ is a connected $\alpha$-critical graph, if $V_1,V_2$ partition
 $V(G)$, and if $\alpha(G)=\alpha(G[V_1])+\alpha(G[V_2])$, then either $V_1$ or $V_2$ is empty.
 \end{prop}
 
 \begin{proof}
 Suppose that both $V_1$ and $V_2$ are nonempty. Since $G$ is connected, there
  must be some pair of adjacent vertices $v\in V_1$ and $w\in V_2$;
  and since $G$ is $\alpha$-critical, it is also the case that
  $\alpha(G-vw)=\alpha(G)+1$. But 
  $\alpha(G-vw)\leq\alpha(G[V_1])+\alpha(G[V_2])$. So 
  $\alpha(G)+1=\alpha(G-vw)\leq\alpha(G[V_1])+\alpha(G[V_2])=\alpha(G)$,
  which is impossible. 
\end{proof}
 
 \begin{cor}\label{cor:XXc}
 If $G$ is a connected $\alpha$-critical graph with maximum critical
 independent set $I$, and $X=I\cup N(I)$, then either $X$ or $X^c$ is
 empty. 
 \end{cor}
 
 \begin{proof}
 The sets $X$ and $X^c$ partition $V(G)$, and Theorem~\ref{idt} yields
 $\alpha(G)=\alpha(G[X])+\alpha(G[X^c])$. 
 \end{proof}
 
 \begin{cor}\label{cor:either}
 If a connected graph $G$ is $\alpha$-critical, then $G$ is either KE or 2-bicritical.
 \end{cor}

\subsection*{Extension to unmatchable graphs}

We have so far focused on graphs with perfect matchings. While random
graphs of even order indeed have perfect matchings
(see, e.g., \cite[p.\,178]{Boll01}), it's natural to ask 
what can be said about graphs of odd order or other
graphs admitting no perfect matching. Do our results yield anything
about the independence structure of unmatchable graphs? Indeed they
do. 

Deming~\cite{Demi79} gave a useful extension of a general graph
$G$ to one with a perfect matching. Let $M$ be a maximum matching 
in $G$, and let $S$ be the set of $M$-unsaturated vertices of $G$. 
For each  $v\in S$, add a new
vertex $v'$ adjacent to $v$, and for each edge $vw$ in $G$, add an
edge $v'w$; denote the resulting graph by $G'$. 
(Notice that the closed neighborhoods $N[v]$ and $N[v']$ in $G'$
are the same for each $v\in S$.)
Now $G'$ has a perfect matching $M'$ (consisting of $M$
together with the edges $vv'$, for $v\in S$), and any maximum
independent set in $G$ is necessarily a maximum independent set in
$G'$ (see the following paragraph). We call the supergraph $G'$ the 
\textit{Deming extension} of $G$ with respect to $M$
and $M'$ its \textit{standard} perfect matching. 

Without specific reference to Deming extensions, two vertices $v$, 
$v'$ in a graph $G$ are \textit{twins} if they have the same closed
neighborhoods; i.e., $N_G[v]=N_G[v']$---the definition implies that $v$
and $v'$ are adjacent. Twins arise in independence theory because a
twin can be removed from $G$ without changing its independence
number; to wit, if $v$ and $v'$ are twins in $G$ and $\widetilde{G}=G-v'$,
then $\alpha(\widetilde{G})=\alpha(G)$ (for if $I$ is a maximum independent set of
$G$ containing $v'$, then $v\not\in I$ and $I'=I-v'+v$ has $|I'|=|I|$
and is a maximum independent set of $\widetilde{G}$). Because a Deming
extension $G'$ is obtained from $G$ by adding certain twins, this
proves the following

\begin{prop}\label{prop:extension}
If $G'$ is a Deming extension of $G$ (with respect to \emph{any} maximum
matching $M$), then $\alpha(G')=\alpha(G)$. 
\end{prop}

\noindent
Theorem~\ref{thm:twins} below expands on this observation. Using 
Proposition~\ref{prop:extension}, it is easy to check that $G'$ is KE
if and only if $G$ is KE.

So the theory developed in this paper for matchable graphs can be
applied to a Deming extension $G'$: $G'$ can be decomposed into Deming
subgraphs and a KE subgraph; each of these parts defines a
corresponding part in $G$---that will not necessarily be a Deming
subgraph or a KE subgraph of $G$, but will be a subgraph of one of
these in $G'$.  These subgraphs interact with their parent graphs in
convenient ways as the next two results show.

\begin{prop}
\label{alpha_same_in_G_and_G'}
Each Deming subgraph $D'$ in a Deming extension $G'$ of a graph $G$ and corresponding
subgraph $D$ of $D'$ in $G$ satisfy $\alpha(D')=\alpha(D)$.  
\end{prop}

\begin{proof}
Let $D'$ be in $G'$, and let $D=D'\cap G$ (that is, the part of $D'$ that's in $G$
or, more formally, $D=G[V(D')\cap V(G)]$, the graph induced on the
common vertices). If $D=D'$ then there is nothing to show. 
So assume that there is a vertex 
$v'\in V(D')\smallsetminus V(D)$ and also that $v'$ is in
some maximum independent set $I'$ of $D'$. 
The edge $vv'$ is by definition in the standard matching $M'$,
and, by the design of Deming's Algorithm,
$vv'$ is an edge in $D'$. Also, $v$ cannot be in $I'$ because $v'$ is
in this independent set.
Now $v$ and $v'$ have the same neighbors in
$G'$ and thus the same neighbors in $D'$. So $v'$ can be replaced with
$v$ to create a new maximum independent set $I=I'-v'+v$. As this can
be done for every such vertex, we obtain a maximum independent set in
$D'$ consisting only of vertices in $D$. And we can conclude that
$\alpha(D)=\alpha(D')$. 
\end{proof}

An analogous argument gives a parallel result for the KE-part of $G'$.

\begin{cor}
\label{alpha_same_in_R_and_R'}
The KE subgraph $R'$ in a Deming extension $G'$ of a graph $G$ and corresponding
subgraph $R$ of $R'$ in $G$ satisfy $\alpha(R')=\alpha(R)$.  
\end{cor}

Thus, the independence structure of a general graph $G$ with no perfect
matching is mirrored (in the sense given) in the independence
structure of the matchable supergraph $G'$. 

Moving forward, the reader should keep in mind the compendium of
perfect matchings introduced by now. We have the standard perfect
matching of $G'$, which may or may not coincide with an `initial'
perfect matching used to initialize the computation of a Deming
decomposition. It's possible that either of these matchings may not
induce perfect matchings for the parts of the Deming
decomposition. But these parts each have perfect matchings, the union
of which yields a (third) perfect matching for the parent graph. 
In considering Deming decompositions, it is usually most convenient to
work with these last `Deming decomposition induced' perfect matchings.

Proposition~\ref{prop:extension} and the results following it show
that Deming extensions $G'$---which are motivated by matching
considerations---preserve simple independence properties of a given
graph $G$. The next result implies that the connection runs still
deeper, to $\alpha$-criticality. Thus, we believe that Deming
decompositions of $\alpha$-critical graphs will be applicable to
further investigation of this important graph class.

The result bears a resemblance to the `Replication Lemma' of 
Lov\'{a}sz and Plummer~\cite[Lemma~12.2.2]{LovaPlum86}, which
substitutes the property of being a perfect graph for being
$\alpha$-critical. (Of course, graphs can have either of these
properties without the other; consider a $4$-cycle [not
$\alpha$-critical] and a $5$-cycle [not perfect].  So the theorems are
not related.)  This result follows from a theorem of
Wessel~\cite{Wess70}, but we give an independent proof.

\begin{thm}\label{thm:twins}
If $v$ and $v'$ are twin vertices in a graph $G$, and $\widetilde{G}=G-v'$,
then $\widetilde{G}$ is $\alpha$-critical if and only if $G$ is $\alpha$-critical. 
\end{thm}

\begin{proof}
We use a couple of times the well-known fact that
if a graph $H$ is $\alpha$-critical and $xy\in E(H)$, then there is a
maximum independent set of $H$ containing $x$ and every maximum
independent set of $H-xy$ contains both $x$ and $y$ (\cite[Problem~8.12]{Lova79a}). 

Suppose first that $\widetilde{G}$ is $\alpha$-critical. Let $xy$ be an edge in
$G$. If $xy$ is an edge in $\widetilde{G}$, then
$\alpha(\widetilde{G}-xy)=\alpha(\widetilde{G})+1=\alpha(G)+1$. It remains to show that
$\alpha(G-xy)=\alpha(\widetilde{G}-xy)$. Since $\widetilde{G}-xy$ is an induced subgraph of
$G-xy$, we have $\alpha(\widetilde{G}-xy)\leq\alpha(G-xy)$. Now let $I$ be a
maximum independent set in $G-xy$. If $v'$ is not in $I$, then $I$ is
an independent set in $\widetilde{G}-xy$ and 
$\alpha(G-xy)=|I|\leq\alpha(\widetilde{G}-xy)$. If $I$ does contain $v'$, then
$I'=I-v'+v$ must also be 
an independent set in $\widetilde{G}-xy$ with cardinality $|I|$, and
$\alpha(G-xy)=|I'|\leq \alpha(\widetilde{G}-xy)$. 

If $xy$ is not an edge in $\widetilde{G}$, then either $x$ or $y$ must be
$v'$. Suppose that $x=v'$. Then $y$ is either $v$ or a common neighbor of
$v$ and $v'$.  Suppose that $y=v$. Since $\widetilde{G}$ is $\alpha$-critical,
there is a maximum independent set $I'$ of $\widetilde{G}$ containing $v$. Then
$I=I'+v'$ is an independent set in $G-v'v=G-xy$. So
$\alpha(G)+1=\alpha(\widetilde{G})+1=|I'|+1=|I|\leq\alpha(G-xy)$. Suppose then
that $y$ is a common neighbor of $v$ and $v'$. Then $vy$ is an
edge in $\widetilde{G}$, and since $\widetilde{G}$ is $\alpha$-critical, we have
$\alpha(\widetilde{G}-vy)=\alpha(\widetilde{G})+1=\alpha(G)+1$. So it's enough to show that
\begin{equation}
\label{key-for-GG-prime-crit}
\alpha(G-v'y)=\alpha(\widetilde{G}-vy).
\end{equation}
Let $I$ be a maximum independent set in
$G-v'y$. If $I$ does not contain $v'$, then $I$ is independent in $\widetilde{G}$
and hence independent in $\widetilde{G}-vy$. So 
$\alpha(G-v'y)=|I|\leq\alpha(\widetilde{G}-vy)$. If $I$ does contain $v'$, it cannot contain any
neighbor of $v'$ except $y$; then $I'=I-v'+v$ must be an independent
set in $\widetilde{G}-vy$. Thus $\alpha(G-v'y)=|I|=|I'|\leq\alpha(\widetilde{G}-vy)$.
We've shown that in either case---i.e.\ $I$ does not, or does, contain
$v'$---the relation
\begin{equation}
\label{half-key-for-GG-prime-crit}
\alpha(G-v'y)\leq\alpha(\widetilde{G}-vy)
\end{equation}
holds. To establish the reverse inequality, consider a maximum independent
set $I$ in $\widetilde{G}-vy$. As $\widetilde{G}$ is $\alpha$-critical, the set $I$ contains
both of $v$, $y$. Thus, $I'=I-v+v'$ must be independent
in $G-v'y$. So now we have $\alpha(G-v'y)\geq |I'|=|I|=\alpha(\widetilde{G}-vy)$,
and combining this with (\ref{half-key-for-GG-prime-crit})
gives (\ref{key-for-GG-prime-crit}).

Conversely, suppose that $G$ is $\alpha$-critical. Let $xy$ be an edge
in $\widetilde{G}$. Since $\widetilde{G}$ is a subgraph of $G$, the edge $xy$ lies in $G$ and
$\alpha(G-xy)=\alpha(G)+1$. As $\alpha(G)=\alpha(\widetilde{G})$, it's enough
to show that $\alpha(G-xy)=\alpha(\widetilde{G}-xy)$. But $\widetilde{G}-xy$ is an induced subgraph
of $G-xy$, whence $\alpha(\widetilde{G}-xy)\leq\alpha(G-xy)$. Now let $I$ be a maximum
independent set in $G-xy$, so that $|I|=\alpha(G-xy)$.  If $I$ does not
contain $v'$, then $I$ is an independent set in $\widetilde{G}-xy$ and
$\alpha(\widetilde{G}-xy)\geq |I|=\alpha(G-xy)$. If $I$ does contain $v'$, then
$I'=I-v'+v$ must also be an independent set in $\widetilde{G}-xy$ with
cardinality $|I|$, and $\alpha(\widetilde{G}-xy)\geq |I'|=\alpha(G-xy)$. 
\end{proof}

It is a curiosity that a large number of $\alpha$-critical graphs
appearing in the literature have twin vertices; these can be reduced
by Theorem~\ref{thm:twins} to smaller twin-free $\alpha$-critical
graphs. For instance, all complete graphs are in this sense
`twin-equivalent'. $K_1$ is $\alpha$-critical, and since $K_2$
is the extension of $K_1$ by a twin, it is $\alpha$-critical, and so on:
since $K_{n-1}$ is $\alpha$-critical and $K_n$ is a twin-extension of
$K_{n-1}$, it too is $\alpha$-critical. 

\begin{cor}\label{thm:criticalextension}
If $G'$ is a Deming extension of a graph $G$, then
$G'$ is $\alpha$-critical if and only if $G$ is $\alpha$-critical. 
\end{cor}

\begin{proof}
A Deming extension can be viewed as the successive addition of a twin
vertex $v'$ for each vertex $v$ that is left unsaturated by a
maximum matching. 
\end{proof}

\subsection*{Gallai class number}

The \textit{Gallai class number} $\delta=n-2\alpha$---a
measure of complexity of critical graphs---has been
central in $\alpha$-critical graph investigations; see, e.g., \cite[Chapter~12]{LovaPlum86}.
Chv\'{a}tal conjectured, and Sewell and Trotter proved~\cite{SeweTrot93},
that every $\alpha$-critical graph with $\delta\geq 2$
contains an even subdivision of $K_4$.  
Matchable graphs $G$ satisfy $n=2\nu$, so we have
$\delta=2\nu-2\alpha$. Since $\nu\geq\alpha$ for 
such graphs, we have $\delta\geq 0$ in this case,
and $\delta$ is necessarily even. Also,
$\delta=0$ if and only if  $G$ is KE, and $\delta=2$ if and only if 
$\alpha=\nu-1$. 

\begin{prop}\label{prop:criticalKE}
If $G$ is a connected $\alpha$-critical KE graph with a perfect
matching $M$, then $G$ is $K_2$. 
\end{prop}

\begin{proof}
The hypotheses give $\alpha=\nu$. For a contradiction, suppose that
$\nu>1$, so that $|M|\geq 2$ and, since $G$ is connected, it contains a
non-$M$ edge $e$. Now criticality gives
$\alpha(G-e)=\alpha(G)+1=\nu(G)+1=\nu(G-e)+1$, which is
impossible. Thus $\nu=1$, and since $G$ has a perfect matching,
$G=K_2$. 
\end{proof}

\begin{cor}\label{cor:K2}
If $G$ is a connected $\alpha$-critical KE graph, then $G$ is $K_1$ or $K_2$.
\end{cor}

\begin{proof}
Let $G$ be a nontrivial connected $\alpha$-critical KE graph. 
If $G$ is matchable, then Proposition~\ref{prop:criticalKE} shows that
$G=K_2$, so we may assume that $G$ contains no perfect matching. We
shall see that this leads to a contradiction.
If we form the Deming extension $G'$ of $G$, then
Corollary~\ref{thm:criticalextension} shows that $G'$ is
$\alpha$-critical and Proposition~\ref{prop:extension}
that $\alpha(G')=\alpha(G)$. By the definition of a Deming extension,
in forming $G'$, we add $(n(G)-2\nu(G))$ new vertices to $G$; so
$n(G')=n(G)+(n(G)-2\nu(G))$ and $\nu(G')=\nu(G)+(n(G)-2\nu(G))$.
With $G$ being KE, we have $\alpha(G)+\nu(G)=n(G)$, and it follows that
$\alpha(G')+\nu(G')=n(G')$, so that $G'$ is also KE. But then 
Proposition~\ref{prop:criticalKE} implies that $G'$ is $K_2$, which
contradicts the fact that $G$ is nontrivial.
\end{proof}

\begin{cor}\label{cor:KE2bi}
If $G$ is a connected $\alpha$-critical graph, then $G$ is either
$K_1$, $K_2$, or is $2$-bicritical. 
\end{cor}

\begin{proof}
Corollary~\ref{cor:either} shows that either $G$ is KE or $G$ is
$2$-bicritical. Now Corollary~\ref{cor:K2} shows that if $G$ is KE, then it
is $K_1$ or $K_2$.  
\end{proof}

We now obtain new proofs of three basic facts about $\alpha$-critical
graphs. For the first, cf.~\cite[Problem~8.21]{Lova79a}.

\begin{cor}
\label{cor:I-nbr-lb}
If $I$ is an independent set in an $\alpha$-critical graph $G$ with
no isolated vertices, then $|N(I)|\geq |I|$.
\end{cor}

\begin{proof}
It's enough to establish the conclusion for a component $H$ of $G$ and 
an independent set $I_H=I\cap V(H)$ of $H$. 
We may assume that $I_H$ is nonempty. 
As $H$ is also $\alpha$-critical, the hypotheses and
Corollaries~\ref{cor:either}, \ref{cor:K2} together imply that
$H$ is either $K_2$ or $2$-bicritical. If $H$ is $K_2$, then $I_H$
consists of a single vertex which has a neighbor. If $H$ is
$2$-bicritical, then $|N(I_H)|>|I_H|$ is part of the defining
condition for such graphs. 
\end{proof}

For the second, cf.~\cite[Corollary~12.1.11]{LovaPlum86}.

\begin{cor}
Every $\alpha$-critical graph $G$ without isolated vertices contains a
perfect $2$-matching.
\end{cor}

\begin{proof}
As in the proof of Corollary~\ref{cor:I-nbr-lb}, each component of
$G$ is either $K_2$ or $2$-bicritical. Both cases admit perfect
$2$-matchings, all of which can be assembled into a perfect
$2$-matching of $G$ (\cite[Corollary~6.2.2]{LovaPlum86} addresses the
$2$-bicritical case).
\end{proof}

For the third, cf.~\cite[Corollary~12.1.12]{LovaPlum86}.

\begin{cor}
If $G$ is $\alpha$-critical and contains no isolated vertices, then
$\alpha(G)\leq n/2$.
\end{cor}

\begin{proof}
Since $\alpha$ is additive across $G$'s
components, it again suffices to obtain the bound for each
component $H$. If $H=K_2$, then this is clear, and if $H$ is
$2$-bicritical, then the bound follows from the defining condition 
($|N_H(I)|>|I|$) applied to a maximum independent set $I$ of $H$.
\end{proof}

The following result is indicated in \cite[p.\,453]{LovaPlum86} and
obtained using a result of Hajnal~\cite{Hajn65} linking the Gallai class
number $\delta$ to maximum degree. We include it because it also
follows from the theory developed here.

\begin{prop}
\label{prop:Andrasfai-pre-helper}
If a graph $G$ is connected, $\alpha$-critical, and $\delta(G)=0$,
then $G=K_2$.   
\end{prop}

\begin{proof}
As $0=\delta(G)=n-2\alpha$ and $\alpha\geq 1$, we see that $G$ has 
even order $n\geq 2$.
As $G$ is connected and $\alpha$-critical, Corollary~\ref{cor:either}
implies that $G$ is either KE or $2$-bicritical. 
But $G$ cannot be $2$-bicritical: if $I$ is a maximum independent set
(with $\alpha=|I|$), it is impossible for both $n=2\alpha$ and
$|N(I)|>|I|$.  Then Corollary~\ref{cor:K2} implies that $G$ is $K_1$
or $K_2$; and since $n\geq 2$, we have $G=K_2$. 
\end{proof}

The following theorem about graphs with $\alpha=\nu-1$ implies a
characterization of those connected $\alpha$-critical graphs with
Gallai class number $2$ (Theorem~\ref{cor:delta2} below).

\begin{thm}\label{thm:numinus1}
A graph $G$ with a perfect matching $M$ satisfies
$\alpha(G)=\nu(G)-1$ if and only if every Deming decomposition of $G$ 
consists of a KE graph $R$ and a single Deming subgraph $D$ (with perfect
matching $M_D=M\cap E(D)$) such that for some edge $xy\in M_D$, the
graph $G-\{x,y\}$ is KE. 
\end{thm}

\begin{proof}
Write a given Deming decomposition (with respect to $M$) as
$\{\{ B_i\}_{i=1}^r, \{K_j\}_{j=1}^{\ell},R\}$. 
Theorem~\ref{thm:dd} shows that 
$\alpha(G)\leq\nu(G)-(r+\ell)$.   

Suppose first that $\alpha(G)=\nu(G)-1$. Then $r+\ell= 1$, so either $r=1$
or $\ell=1$, whence $G$ has exactly one Deming subgraph $D$, and it must be the
case that $\alpha(D)=\nu(D)-1$ (see Theorem~\ref{thm:dd}, parts~(1,2,5)).
Let $I$ be a maximum independent set of $G$ (so $|I|=\nu(G)-1$), and 
let  $I_D=I\cap V(D)$ and $I_R=I\cap V(R)$; $I_D$ is then a maximum
independent set of $D$ and $I_R$ is a maximum independent set of $R$.
Note that no more than one vertex of $I_D$ is incident to any edge
of $M_D$, and---because $\alpha(D)=\nu(D)-1$---there must be one edge
$xy\in M_D$ that is not incident with any vertex in $I_D$. Let
$G'=G-\{x,y\}$. By construction, 
$I\subseteq V(G')$, and $M'=M-\{xy\}$ is a perfect matching of $G'$ with
cardinality $|I|$. Thus $G'$ is KE. 

Now suppose that every Deming decomposition consists
of a KE graph $R$ and a single Deming subgraph $D$ (with perfect
matching $M_D=M\cap E(D)$) such that for some edge $xy\in M_D$, the graph
$G'=G-\{x,y\}$ is KE. 
Since $M$ is a perfect matching of $G$ and $xy\in M$, the graph $G'$
is also matchable and $\nu(G')=\nu(G)-1$. Now $G'$, being KE, has
$\alpha(G')=\nu(G')$, which implies that $G'$ has an independent set
$I\subseteq V(G')$ of order $\nu(G')$. Thus, $I$ consists of exactly
one end of each edge in $M-xy$, and therefore $I$ is also independent
in $G$. So we have
\begin{equation}
\label{alpha-lb}
\alpha(G)\geq |I|=\nu(G')=\nu(G)-1.
\end{equation}
But $G$ is not KE because its Deming decomposition contains a non-null
Deming subgraph ($D$). This means that the chain (\ref{alpha-lb}) is sharp;
i.e., $\alpha(G)=\nu(G)-1$.
\end{proof}

We're soon to recover the result of Andr\'asfai~\cite{Andr67} on connected
$\alpha$-critical graphs with $\delta=2$ (to which we alluded prior to
the statement of Theorem~\ref{thm:numinus1}); first we need to show
that such graphs are matchable.

\begin{lem}
\label{lem:Andrasfai-helper}
If a graph $G$ is connected, $\alpha$-critical, and $\delta(G)=2$,
then $G$ contains a perfect matching.
\end{lem}

\begin{proof}
As $2=\delta(G)=n-2\alpha$ and $\alpha\geq 1$, we see that $G$ has 
even order $n\geq 4$. With $G$ being connected and $\alpha$-critical,
Corollary~\ref{cor:either} shows that $G$ is either KE or $2$-bicritical. 
But it's not the first of these, for if so, then
Corollary~\ref{cor:K2} would give $G\in\{K_1,K_2\}$ and thus
contradict $n\geq 4$. 
Since  $G$ is $2$-bicritical, it contains a spanning subgraph $H$
consisting of a matching $M$ together with an even number of odd
cycles $C_1,C_2,\ldots,C_{\ell}$
(cf.\ \cite[Corollary~6.2.2]{LovaPlum86}, where $H$ is termed a `basic
perfect $2$-matching' of $G$). 
As $n\geq 4$ and $G$ is connected, $G$ must contain an
edge $e\not\in E(H)$. Because $G$ is $\alpha$-critical, we have
\[
\alpha(G-e) = \alpha(G)+1 = \frac{n-2}{2}+1 = \frac{n}{2}.
\]
Now consider an independent set $I$ in $G-e$ with $|I|=n/2$. Since 
$e\not\in E(H)$, the subgraph $H$ spans $G-e$, and so the set $I$ is also
independent in $H$. Such a set can occupy at most one end of each edge
of $M$ and at most $(|C_i|-1)/2$ vertices of each cycle $C_i$. Thus,
\begin{equation}
\label{I-upper-bnd}
|I|\leq |M| + \sum_{i=1}^{\ell}\frac{|C_i|-1}{2} = \frac{n}{2}-\frac{\ell}{2},
\end{equation}
and since $|I|=n/2$, the bound (\ref{I-upper-bnd}) shows that 
$\ell=0$. That is, $H$ contains no cycle
components and is therefore a perfect matching of $G$. 
\end{proof}

Perhaps surprisingly, the hypotheses of
Lemma~\ref{lem:Andrasfai-helper} guarantee not just a perfect matching
in $G$ but an exact, concise description of this graph. At this point,
it may be instructive to review our earlier and much simpler
Proposition~\ref{perfect}.

\begin{thm}[Andr\'asfai, \cite{Andr67}]\label{cor:delta2}
If a graph $G$ is connected, $\alpha$-critical, and $\delta(G)=2$, then $G$
is an even subdivision of $K_4$.
\end{thm}

\begin{proof}
Denote by $M$ the perfect matching guaranteed by Lemma~\ref{lem:Andrasfai-helper}.
We showed above that $\delta(G)=2$ is equivalent to $\alpha(G)=\nu(G)-1$. 
Then Theorem~\ref{thm:numinus1}
implies that  a  Deming decomposition of $G$ consists of a KE graph
$R$ and a single Deming subgraph $D$ (with perfect matching 
$M_D=M\cap E(D)$). 

While $R$ (which is KE) may not have any vertices, $D$ must (if 
$V(D)=\varnothing$, then $\alpha(G)=\alpha(R)=\nu(R)=\nu(G)$, contrary to
the fact that $\alpha(G)=\nu(G)-1$). Assume then that
$V(R)\not=\varnothing$. The relation 
$\alpha(G)=\nu(G)-1$ implies that $\alpha(R)=\nu(R)$ and
$\alpha(D)=\nu(D)-1$. And $\nu(G)=\nu(R)+\nu(D)$ implies that
$\alpha(G)=\alpha(R)+\alpha(D)$. Since $G$ is connected, there must be
an edge from some vertex $x$ in $R$ to some vertex $y$ in $D$. 
Notice that $\alpha(G-xy)=\alpha(G)$. (Otherwise,
$\alpha(G-xy)>\alpha(G) =\alpha(R)+\alpha(D)$; but then a maximum 
independent set in $G-xy$  would, by the
Pigeonhole Principle, contain either more than $\alpha(R)$ vertices of
$R$ or more than $\alpha(D)$ vertices of $D$ and so could not be
independent.)
But this contradicts $G$ being $\alpha$-critical. So it must be that
$R$ is a null graph and $G=D$.  

With $D$ being an $\alpha$-critical Deming graph,
Corollary~\ref{cor:alphacriticaldeming} shows that $D(=G)$ is an even
subdivision of $K_4$.
\end{proof}

The following theorem suggests a tantalizing connection between the structure of $\alpha$-critical graphs and Deming decompositions. Is it true that, for instance, any connected $\alpha$-critical graph with a perfect matching has an associated Deming decomposition with a Deming-$K_4$ graph (and hence the required even $K_4$-subdivision)?

\begin{thm}[Sewell and Trotter, \cite{SeweTrot93}]\label{thm:sewell}
If a graph $G$ is connected, $\alpha$-critical, and $\delta(G)\geq 2$, then $G$
contains an even subdivision of $K_4$.
\end{thm}

What more can be said about the independence structure of the
$2$-bicritical subgraph of an arbitrary graph? We assume that it
has  a perfect matching. 
It is a simple consequence of this 
and the defining condition of 2-bicritical graphs that any independent
set can be matched and thus that $\alpha<\nu$. We also assume that
some Deming decomposition of the graph has $k$ Deming graphs (either
Deming-BP or Deming-$K_4$ graphs) $D_i$. We know (from
Theorem~\ref{thm:deming}) 
that these Deming graphs
satisfy $\alpha(D_i)=\nu(D_i)-1$. An interesting feature of the
following claims and related algorithms is that a perfect matching $M$
must be computed only once---it can be reused for all further computations.

The following result generalizes Theorem~\ref{thm:numinus1}---needed
for our investigation of $\alpha$-critical graphs---to matchable graphs
with $k$ Deming subgraphs in a Deming decomposition. 

\begin{thm}
Let $G$ be a graph with a perfect matching $M$ and a Deming
decomposition consisting of a KE graph $R$ and Deming subgraphs
$\{D_i\}_{i=1}^{k}$ with $k\geq 1$. Let $M_i$ be the edges of $M$ in
$D_i$. Then $\alpha(G)=\nu(G)-k$ if and only if there is an edge 
$x_iy_i\in M_i$, for each $i\in [k]$, such that
$G-\{x_1,\ldots,x_k,y_1,\ldots,y_k\}$ is KE. 
\end{thm} 

\begin{proof}
The independence number of $G$ is no more than the sum of the
independence numbers of the KE and Deming subgraphs from a Deming
decomposition. Thus,
\begin{equation}
\label{alpha-ub-gen}
\alpha(G)\leq \alpha(R)+\sum_{i=1}^k\alpha(D_i)=
\nu(R)+\sum_{i=1}^k(\nu(D_i)-1)=\nu(G)-k,
\end{equation}
the first equality following from Theorem~\ref{thm:deming} and the
fact that $R$ is KE.

Suppose first that $\alpha(G)=\nu(G)-k$, and 
let $I$ be a maximum independent set of $G$ (so $|I|=\nu(G)-k$).
With $I_i=I\cap V(D_i)$ and $I_R=I\cap V(R)$, the sets $I_i$ are
maximum independent sets of their respective $D_i$'s, and $I_R$ is a
maximum independent set of $R$.  
Note that no more than one vertex of $I_i$ is incident to any
edge in $M_i$, and---because $\alpha(D_i)=\nu(D_i)-1$---there 
must be one edge $x_iy_i\in M_i$
that is not incident with any vertex in $I_i$. Let
$G'=G-\{x_1,\ldots,x_k,y_1,\ldots,y_k\}$. By construction, 
$I\subseteq V(G')$ and $M'=M-\{x_1y_1,\ldots,x_ky_k\}$ is a perfect
matching of $G'$ with cardinality $|I|$. Thus $G'$ is KE. 

Suppose now that there is an edge $x_iy_i\in M_i$ for each $i\in [k]$
such that $G'=G-\{x_1,\ldots,x_k,y_1,\ldots,y_k\}$ is KE.  Since $G'$
has the perfect matching $M'=M-\{x_1y_1,\ldots,x_ky_k\}$ with
$|M'|=|M|-k=\nu(G)-k$, it follows that
$\alpha(G')=\nu(G')=\nu(G)-k$, which shows that 
$G'$ has an independent set
$I\subseteq V(G')$ of order $\nu(G)-k$. Thus, $I$ consists of exactly
one end of each edge in $M-\{x_1y_1,\ldots,x_ky_k\}$, and so $I$ is 
also independent in $G$. We therefore have
\begin{equation}
\label{alpha-lb-gen}
\alpha(G)\geq |I|=\nu(G)-k,
\end{equation}
and combining (\ref{alpha-ub-gen}) with (\ref{alpha-lb-gen}) gives $\alpha(G)=\nu(G)-k$.
\end{proof}

It may not be the case that $\alpha(G)=\nu(G)-1$ for a graph $G$ with
a perfect matching $M$ and a Deming decomposition with a single Deming
subgraph $D$. Theorem~\ref{thm:numinus1} identifies a condition that
guarantees this; namely, for some edge $xy\in M_D$, the graph $G-\{x,y\}$ is KE. 
Failing this, $\alpha(G)=\nu(G)-k$ for some $k>1$. Can this
$k$ be efficiently identified? This (open) question provides a natural
segue to our closing section. 
  
\section{(More) Open Problems}
\label{sec:open}

The motivating questions of this investigation were: 
\begin{enumerate}
\item If a graph is Egerv\'{a}ry, is there an efficiently checkable certificate?
\item Is there an NP $\cap$ co-NP description (i.e., a good
characterization) of graphs that have no disjoint pair of odd
cycles?  Or is it NP-complete to decide? 
\item Can Egerv\'{a}ry graphs be recognized in polynomial-time?
\item Is it possible to identify a maximum independent set in an Egerv\'{a}ry graph in polynomial-time?
\end{enumerate} 
Our investigations herein raise a new question:
\begin{enumerate}
\addtocounter{enumi}{4}
\item If $G$ is a graph with a perfect matching $M$ and a Deming
 decomposition $\{\{ K_j\}_{j=1}^{\ell},R\}$ (with respect to $M$), is
 it true that $\alpha(G)= \sum_{j=1}^{\ell} \alpha(K_j)+ \alpha(R)$? 
\end{enumerate}

\noindent
These questions remain open.

\subsection*{Acknowledgements}
Heartfelt thanks to Jack Edmonds for inspiring us in the conception,
sweat, and production of this paper. Its seeds were planted at the
2012 SIAM Conference on Discrete Mathematics in Halifax, Canada, when
the three of us participated in a mini-symposium on KE graphs. Out of
this grew a joint project between we three; see
\cite{Edmonds-ACCOTA}. Beyond the inspiration and collaboration, 
Jack coined the term `Egerv\'{a}ry', wrote the first draft of
Section\,\ref{sec:constructing},  and shared many other related
results and wisdom.

\bibliographystyle{plain}
\bibliography{kayll-larson}

\def\cprime{$'$} \def\cprime{$'$} \def\cprime{$'$} \def\cprime{$'$}
\begin{thebibliography}{10}

\bibitem{Andr67}
B.~Andr\'{a}sfai.
\newblock On critical graphs.
\newblock In {\em Theory of Graphs, International Symposium, Rome, July 1966},
  pages 9--19. Gordon and Breach, New York, 1967.

\bibitem{BeinHaraPlum67}
L.~Beineke, F.~Harary, and M.~Plummer.
\newblock On the critical lines of a graph.
\newblock {\em Pacific Journal of Mathematics}, 22(2):205--212, 1967.

\bibitem{Boll01}
B.~Bollob{\'a}s.
\newblock {\em Random graphs}, volume~73 of {\em Cambridge Studies in Advanced
  Mathematics}.
\newblock Cambridge University Press, Cambridge, second edition, 2001.

\bibitem{CarvKothWangLin20}
M.~H. de~Carvalho, N.~Kothari, X.~Wang, and Y.~Lin.
\newblock Birkhoff--von {N}eumann graphs that are {PM}-compact.
\newblock {\em SIAM J. Discrete Math.}, 34(3):1769--1790, 2020.

\bibitem{Demi79}
R.~W. Deming.
\newblock Independence numbers of graphs---an extension of the
  {K}oenig-{E}gervary theorem.
\newblock {\em Discrete Math.}, 27(1):23--33, 1979.

\bibitem{Edmo65}
J.~Edmonds.
\newblock Paths, trees, and flowers.
\newblock {\em Canad. J. Math.}, 17:449--467, 1965.

\bibitem{Edmonds-ACCOTA}
J.~Edmonds, M.~Kayll, and C.~Larson.
\newblock Fractional vertex packing and covering.
\newblock ACCOTA 2012 (International Workshop on Combinatorial and
  Computational Aspects of Optimization, Topology and Algebra), 3~December
  2012.
\newblock Huatulco, Oaxaca, M\'{e}xico.

\bibitem{ErdoGall61}
P.~Erdos and T.~Gallai.
\newblock On the minimal number of vertices representing the edges of a graph.
\newblock {\em Publ. Math. Inst. Hungar. Acad. Sci}, 6(18):1--203, 1961.

\bibitem{Hajn65}
A.~Hajnal.
\newblock A theorem on {$k$}-saturated graphs.
\newblock {\em Canadian J. Math.}, 17:720--724, 1965.

\bibitem{KawaOzek13}
K.~Kawarabayashi and K.~Ozeki.
\newblock A simpler proof for the two disjoint odd cycles theorem.
\newblock {\em J. Combin. Theory Ser. B}, 103(3):313--319, 2013.

\bibitem{Kayl10}
P.~M. Kayll.
\newblock K\"{o}nig-{E}gerv\'{a}ry graphs are non-{E}dmonds.
\newblock {\em Graphs and Combinatorics}, 26(5):721--726, 2010.

\bibitem{Lars11}
C.~E. Larson.
\newblock The critical independence number and an independence decomposition.
\newblock {\em European J. Combin.}, 32(2):294--300, 2011.

\bibitem{Lova79a}
L.~Lov{\'a}sz.
\newblock {\em Combinatorial problems and exercises}.
\newblock North-Holland Publishing Co., Amsterdam, 1979.

\bibitem{LovaPlum86}
L.~Lov{\'a}sz and M.~D. Plummer.
\newblock {\em Matching {T}heory}, volume 121 of {\em North-Holland Mathematics
  Studies}.
\newblock North-Holland Publishing Co., Amsterdam, 1986.
\newblock Annals of Discrete Mathematics, 29.

\bibitem{Pull79}
W.~R. Pulleyblank.
\newblock Minimum node covers and 2-bicritical graphs.
\newblock {\em Mathematical programming}, 17(1):91--103, 1979.

\bibitem{SeweTrot93}
E.~C. Sewell and L.~E. Trotter.
\newblock Stability critical graphs and even subdivisions of k4.
\newblock {\em Journal of Combinatorial Theory, Series B}, 59(1):74--84, 1993.

\bibitem{Slil07}
D.~Slilaty.
\newblock Projective-planar signed graphs and tangled signed graphs.
\newblock {\em J. Combin. Theory Ser. B}, 97(5):693--717, 2007.

\bibitem{Ster79}
F.~Sterboul.
\newblock A characterization of the graphs in which the transversal number
  equals the matching number.
\newblock {\em J. Combin. Theory Ser. B}, 27(2):228--229, 1979.

\bibitem{Wess70}
W.~Wessel.
\newblock On the problem of determining whether a given graph is edgecritical
  or not.
\newblock In {\em Combinatorial theory and its applications, {III} ({P}roc.
  {C}olloq., {B}alatonf\"{u}red, 1969)}, pages 1123--1139, 1970.

\end{thebibliography}

\end{document}